\documentclass[a4paper,11pt,fleqn]{article}

\usepackage{graphicx}
\usepackage{relsize}
\usepackage{amsfonts}
\usepackage{verbatim}
\usepackage{xifthen}
\usepackage{xfrac}
\usepackage{xargs}
\usepackage{xcolor}
\usepackage[hidelinks]{hyperref}
\usepackage{physics}
\usepackage{amsmath} % For advanced math typesetting
\usepackage{amsthm}  % For defining theorem-like environments
\usepackage{amssymb} % For additional mathematical symbols
\usepackage{listings} % For including code with syntax highlighting
\usepackage{tikz}
\usepackage{authblk}
\usepackage{mathtools}
\usepackage[a4paper,top=1in,bottom=1in,left=1in,right=1in]{geometry}
\usepackage{bm}
\usepackage{subcaption}
\usepackage{float}
\usepackage{algorithm}
\usepackage{algpseudocode}
\usepackage{tocloft}
\usepackage{booktabs}

\usetikzlibrary{positioning}
\usetikzlibrary{arrows}
\usetikzlibrary{decorations.pathmorphing}

\newcommand {\magn}[1]{\left| #1 \right|}
\renewcommand {\vec}[1]{\boldsymbol{#1}}

\newcommand \I{\mathrm{i}}
\newcommandx\D[1][1=]{\mathrm{d}#1}

\newcommandx {\equals}[3][1={=}]{\underbrace{#2}_{#1#3}}

\newcommandx {\CouLombV}[5][3=,5=]{\frac{#1}{\magn{\vec{#2}_{#3}-\vec{#4}_{#5}}}}
\newcommandx {\CouLombF}[4][2=,4=]{\frac{\vec{#1}_{#2}-\vec{#3}_{#4}}{\magn{\vec{#1}_{#2}-\vec{#3}_{#4}}^3}}

\newcommandx {\GFun}[4][2=,4=]{G\left(\vec{#1}_{#2},\vec{#3}_{#4}\right)}
\newcommandx {\der}[2][1=]{\frac{\D[#1] }{\D[#2]}}
\newcommandx {\pder}[2][1=]{\frac{\partial #1}{\partial #2}}
\newcommand{\bb}[1]{\left(#1\right)} 
\renewcommand{\sb}[1]{\left[#1\right]}

\renewcommand{\Im}[1]{\mathrm{Im}(#1)}
\renewcommand{\Re}[1]{\mathrm{Re}(#1)}

\renewcommand{\abs}[1]{\lvert #1 \rvert}

\renewcommand{\phi}{\varphi}

\newcommand{\R}{\mathbb{R}}

\theoremstyle{definition}
\theoremstyle{definition} % Style for definitions (upright body)

\theoremstyle{definition} % Style for remarks (upright body, smaller font)
\newtheorem{remark}{Remark}[section]

\begin{document}

\begin{center}
   \begin{minipage}[t]{6.0in}
    We introduce an $O(M)$ algorithm for evaluating the azimuthal Fourier modes $G_{k,m}$\,, $m = 0, \ldots, M$\,, of the three-dimensional Helmholtz Green's function with real wavenumber $k$\,, together with all their first- and second-order derivatives with respect to the cylindrical source and target coordinates. The cost is independent of both the wavenumber and the source-target separation, and high relative accuracy is retained even for modes whose magnitude is exponentially small. The method combines contour deformation at a few boundary modes with a boundary-value formulation of the five-term recurrence in the mode index. Derivative quantities are obtained from stable recurrences, adding only a small constant factor to the cost of $G_{k,m}$ alone. Numerical experiments demonstrate high relative accuracy, linear scaling in $M$, and applications to modal boundary integral equation solvers for axisymmetric acoustic scattering, where the $k$-independent kernel evaluator makes dense per-mode linear algebra the dominant cost.
    \thispagestyle{empty}

  \vspace{ -100.0in}

  \end{minipage}
\end{center}

\vspace{ 2.60in}
\vspace{ 0.50in}

\begin{center}
  \begin{minipage}[t]{4.4in}
    \begin{center}

\textbf{Fast Evaluation of the Azimuthal Fourier Modes of the 3D Helmholtz Green's Function and Their Derivatives}

  \vspace{ 0.50in}

Hanwen Zhang$\mbox{}^{\dagger}$ \\
              \today

    \end{center}
  \vspace{ -100.0in}
  \end{minipage}
\end{center}

\vspace{ 2.00in}

\vfill

\noindent
$\mbox{}^{\dagger}$ Courant Institute of Mathematical Sciences, New York University, New York, NY 10012 \\

\vfill
\eject

\setcounter{tocdepth}{3}
\tableofcontents
\section{Introduction}

In this paper, we introduce an $O(M)$ algorithm for computing the modal Green's functions
\begin{align}
G_{k,m}(r,z,r',z') = \frac{1}{2\pi}	\int_{-\pi}^{\pi} G_k(\abs{\bm{x} - \bm{x}'})\, e^{-\I m\theta}\, \D \theta\,, \quad m = 0, 1, \ldots, M\,,
\label{eq:mgf}
\end{align}
of the three-dimensional Helmholtz Green's function
\begin{align}
	G_k (\abs{\bm{x} -\bm{x}'}) =  \frac{1}{4\pi} \frac{e^{\I k \abs{\bm{x}-\bm{x}'}}}{\abs{\bm{x}-\bm{x}'}}\,,
	\label{eq:helmgf}
\end{align}
where $k$ is the wavenumber and $(r, \phi, z)$, $(r', \phi', z')$ are the cylindrical coordinates of $\bm{x}$ and $\bm{x}'$ with $\theta = \phi - \phi'$\,. The algorithm also produces all first- and second-order derivatives of $G_{k,m}$ with respect to the cylindrical source and target coordinates at little additional cost. The cost is independent of both the wavenumber $k$ and the source-target separation. The algorithm preserves high relative accuracy across all modes, even for those whose magnitude is exponentially small.

The modal Green's functions arise naturally in boundary integral equation (BIE) methods for problems with rotational symmetry \cite{coltonkress2019inverse, kress2014linear}, with applications ranging from acoustic scattering \cite{liu2016efficient} and radar cross-section characterization \cite{youssef1989radar} to plasma physics \cite{oneil2018taylor} and optical metasurface and metalens design \cite{xue2023fullwave, chung2023inverse}. When the boundary $S$ of a domain $\Omega \subset \R^3$ is a surface of revolution, generated by rotating a curve $\gamma$ in the $(r,z)$-plane around the $z$-axis, the integral equation kernel $K$ — built from $G_k$ and its derivatives — depends on the azimuthal coordinates only through $\theta$\,. Expanding in Fourier series and integrating out $\phi'$\,, the Fourier modes decouple into independent one-dimensional integral equations on $\gamma$\,:
\begin{align}
	\sigma_m(r,z) + 2\pi\int_{\gamma} K_m(r,z,r',z')\, \sigma_m(r',z')\, r'\,\D\ell' = f_m(r,z)\,, \quad (r,z) \in \gamma\,,
	\label{eq:bie_modal}
\end{align}
for each mode $m$\,, where $\sigma_m$\,, $K_m$\,, $f_m$ are the $m$-th Fourier coefficients of the density, kernel, and right-hand side, respectively (see, e.g., \cite{epstein2019high, garritano2022efficient}). The modal kernel $K_m$ is a linear combination of $G_{k,m}$ and its first- and second-order derivatives with respect to $(r, z, r', z')$\,. This dimensional reduction is meaningful only when the modal kernel can itself be evaluated efficiently.

Evaluating $G_{k,m}$ poses several numerical challenges. First, the integrand in (\ref{eq:mgf}) oscillates at frequencies set by both $m$ and $k$, so resolving it by naive quadrature requires a number of nodes that grows with both. Second, when $\bm{x}$ and $\bm{x}'$ are close, the integrand becomes near-singular, and standard smooth quadrature fails. Third, for large $m$, $|G_{k,m}|$ becomes exponentially small, so direct quadrature loses relative accuracy through cancellation. Finally, BIE applications often require not only $G_{k,m}$ but also its first- and second-order derivatives, whose integrands carry even stronger near-singular factors.

Several methods have been developed for evaluating modal Helmholtz Green's functions. Contour-deformation methods, introduced in \cite{gustafsson2010accurate} for $G_{k,0}$ and extended to general modes in \cite{garritano2022efficient}, achieve numerically stable evaluation in $O(m)$ time per mode, with cost independent of the wavenumber and source-target geometry; applied mode by mode, however, computing all of $G_{k,0}, \ldots, G_{k,M}$ requires $O(M^2)$ work. An alternative combines the kernel splitting of Helsing and Karlsson \cite{helsing2014kernel} with FFT-based evaluation, as used in \cite{epstein2019high, helsing2017resonances} for axisymmetric Maxwell problems; this evaluates all modes simultaneously, but the per-pair cost grows linearly with the wavenumber. Lai and O'Neil \cite{lai2019fft} extended the kernel-splitting approach to second-order derivatives of the modal Green's function. Matviyenko \cite{matviyenko1995azimuthal} proposed a different approach, expressing $G_{k,m}$ through a five-term recurrence in $m$, but both forward and backward recursions are unstable across a wide range of parameters and the recurrence cannot be used directly. A detailed review of the methods above and related approaches can be found in \cite{garritano2022efficient}. Thus existing approaches do not provide an $O(M)$ evaluation of all modes $G_{k,0}, \ldots, G_{k,M}$ with cost independent of the wavenumber.

Our algorithm combines the advantages of contour deformation and recurrence relations. First, we apply contour deformation only at a small number of boundary modes, inheriting its wavenumber-independence. Second, the remaining modes are recovered by solving the five-term recurrence as a banded boundary-value problem in the mode index, avoiding the instabilities of forward and backward recursion. Third, derivative quantities are obtained from $G_{k,m}$ by stable recurrences, so first- and second-order derivative data require only a constant-factor increase in cost. In the decay regime, where contour deformation no longer provides accurate boundary values, we follow a strategy analogous to Miller's algorithm for evaluating special functions, which retains relative accuracy across many orders of magnitude. In modal BIE solvers for axisymmetric scattering, the resulting $k$-independent kernel evaluator makes dense per-mode linear algebra the dominant cost.

The remainder of the paper is organized as follows. Section~\ref{sec:modalgf} introduces the parametrization of $G_{k,m}$ and identifies three evaluation regimes. Section~\ref{sec:aux} contains the derivative formulas and recurrence relations needed for the algorithm. Section~\ref{sec:contour} summarizes the contour deformation of \cite{garritano2022efficient,gustafsson2010accurate} for a single Fourier mode $G_{k,m}$ and extends it to the first- and second-order derivatives at the same $O(m)$ cost. Section~\ref{sec:allmodes} presents the pentadiagonal boundary-value formulation that yields all $G_{k,m}$ in $O(M)$ time by combining contour deformation with the five-term recurrence, and demonstrates its componentwise accuracy across regimes. Section~\ref{sec:algorithm} assembles the complete algorithm for evaluating $G_{k,m}$\,, $m = 0, \ldots, M$\,, and their first- and second-order derivatives. Section~\ref{sec:numerical} presents numerical experiments: accuracy and timing of the kernel evaluator (Section~\ref{sec:accuracy_timing}), and application to BIE solvers for acoustic scattering from bodies of revolution (Section~\ref{sec:bie}). Section~\ref{sec:conclusion} summarizes the results and discusses future directions.

\section{Modal Green's functions}
\label{sec:modalgf}

This section introduces the dimensionless parametrization used throughout the paper and the three evaluation regimes that determine the numerical strategy. Section~\ref{sec:parametrization} reduces the modal Green's function, via the parametrization of \cite{epstein2019high,garritano2022efficient}, to a form depending only on a frequency parameter $\kappa$ and a geometric parameter $\alpha$\,. Section~\ref{sec:regimes} describes the near-axis, non-decay, and decay regimes.

\subsection{Parametrization and notation}
\label{sec:parametrization}

We begin by expressing the Green's function (\ref{eq:helmgf}) in cylindrical coordinates, so that the modal Green's function (\ref{eq:mgf}) takes the explicit form
\begin{align}
G_{k,m}(r,z,r',z') = \frac{1}{2\pi}	\int_{-\pi}^{\pi} \frac{1}{4\pi} \frac{e^{\I k  \sqrt{r^2 + r'^2 + (z-z')^2 - 2 r r'\cos{\theta}}}}{\sqrt{r^2 + r'^2 + (z-z')^2 - 2 r r'\cos{\theta}}}\cdot e^{-\I m\theta}\, \D \theta\,.
\label{eq:mgf_explicit}
\end{align}
Following \cite{epstein2019high, garritano2022efficient}, we introduce the notation
\begin{align}
	R_0 = \sqrt{r^2 + r'^2 + (z-z')^2} \,,
	\label{eq:r0}
\end{align}
with
\begin{align}
	\kappa = k R_0 \,,   \quad  \mbox{and}\quad \alpha = 2r r'/R_0^2\,.
	\label{eq:kanda}
\end{align}
Here $\kappa$ encodes the oscillatory behavior of the kernel, while $\alpha$ depends purely on the geometry of the source and target locations. It is clear from the definition that $\alpha \in [0,1)$\,. Substituting (\ref{eq:r0}) and (\ref{eq:kanda}) into (\ref{eq:mgf_explicit}) and exploiting the even symmetry of the integrand in $\theta$ (which implies $G_{k,m} = G_{k,-m}$), we obtain
\begin{align}
	G_{k,m}(r,z,r',z') = \frac{1}{4\pi^2} \frac{1}{R_0} \int_{0}^{\pi} \frac{e^{\I\kappa \sqrt{1-\alpha\cos\theta}}}{\sqrt{1-\alpha\cos\theta}} \cdot \cos{ (m\theta)}\,\D \theta\,.
	\label{eq:mgf4}
\end{align}
This symmetry lets us restrict to $m \ge 0$. Closed-form expressions for $G_{k,m}$ in terms of two-variable hypergeometric functions, as well as series expansions in Legendre and Bessel functions, are given in \cite{conway2010exact}; in this paper we work with the integral representation (\ref{eq:mgf4}) throughout. For notational convenience, we write
\begin{align}
	G_m = G_{k,m}\,,
\end{align}
throughout the paper, suppressing the dependence on $k$\,.

It will also be useful to introduce the parameters
\begin{align}
	\beta_- = \sqrt{\frac{{(r-r')^2 + (z-z')^2}}{2rr'}} = \sqrt{1/\alpha-1}\,, \quad\quad  \beta_+ = \sqrt{1/\alpha + 1}\,.
	\label{eq:gmparas}
\end{align}
In the near-singular regime, where $r\to r'$ and $z\to z'$\,, we have $\alpha \to 1^{-}$ and $\beta_- \to 0$\,. In this limit, the integrand in (\ref{eq:mgf4}) develops a weak singularity near $\theta = 0$\,, and $\beta_-$ provides a convenient measure of how close the source and target configurations are to one another. When $\alpha$ is close to $0$, the integrand in (\ref{eq:mgf4}) is smooth with weak $\theta$-dependence. This occurs when the source-target separation is large relative to $\sqrt{rr'}$; for example, when one point is close to the $z$-axis compared with the distance between the source and target.

\subsection{Three evaluation regimes}
\label{sec:regimes}

The evaluation of $G_m$ for $m = 0, 1, \ldots, M$ naturally divides into three regimes, each requiring a different numerical approach. The near-axis regime is characterized by $\alpha$ close to $0$. The remaining two regimes are distinguished by whether the requested mode range extends beyond the transition point at which $|G_m|$ begins to decay exponentially.

\paragraph{Near-axis regime ($\alpha \ll 1$\,, $\kappa\alpha \le 1$).}
In this regime, $G_m$ admits a rapidly convergent power series expansion in $\alpha$\,:
\begin{align}
	G_m = \frac{1}{4\pi R_0}\, e^{\I\kappa} \sum_{n=0}^{N} p_n(\kappa)\, \alpha^{m+2n} + O(\alpha^{m+2N+2})\,,
	\label{eq:nearaxis}
\end{align}
where $p_n$ are polynomials in $\kappa$ that depend on $m$\,. Since $G_m$ decays rapidly with $m$ in this regime, only a small number of modes are significant and the expansion converges quickly. This regime is handled by precomputed Taylor-expansion tables; the remainder of the paper focuses on the non-decay and decay regimes.

\paragraph{Decay regime.}
Outside the near-axis regime, there is a transition point $m^*$ beyond which $|G_m|$ decays exponentially for $m \ge m^*$. A classical stationary phase argument \cite{matviyenko1995azimuthal, garritano2022efficient} gives
\begin{align}
	m^* = \frac{\kappa}{\sqrt{2}}\sqrt{1 - \sqrt{1-\alpha^2}}\,.
	\label{eq:mstar}
\end{align}
For $m>m^*$, evaluating $G_m$ from the integral representation (\ref{eq:mgf4}) loses relative accuracy: the integrand remains of size $O(1)$ and oscillatory, and the exponentially small value of the integral is the result of cancellation among its oscillatory contributions. Since the absolute rounding error incurred by quadrature is governed by the scale of the integrand rather than that of the exponentially small result, relative accuracy is lost even when the oscillation is fully resolved.

\paragraph{Non-decay regime.}
When $M \le m^*$\,, $|G_m|$ does not decay within the range $m = 0, \ldots, M$\,. Therefore, evaluation methods based on the integral representation (\ref{eq:mgf4}) achieve full relative accuracy with sufficiently many quadrature nodes, up to a condition-number factor set by the oscillation of the integrand.

These observations extend unchanged to the derivatives of $G_m$.

\begin{remark}
In this paper, we use $\alpha \le 0.05$ as the threshold for the near-axis regime. For $\alpha$ near $1$, specifically $1 - \alpha \lesssim 10^{-5}$\,, the decay of $|G_m|$ is sufficiently slow that the problem can be treated as non-decay for $M \le 3000$\,, with at most a two-digit loss of relative accuracy.
\end{remark}

\section{Derivatives and recurrence relations}
\label{sec:aux}

This section develops the derivative identities and recurrence relations that make it possible to compute all first- and second-order derivatives of $G_m$ at essentially the same asymptotic cost as $G_m$ itself. Section~\ref{sec:derivatives} introduces the auxiliary parametrization $(a,b)$ following \cite{matviyenko1995azimuthal} and derives integral representations for the first- and second-order derivatives of $G_m$\,, including the regrouped forms that avoid catastrophic cancellation in the near-singular regime. Section~\ref{sec:recurrence} introduces the five-term recurrence for $G_m$ and the associated two-term recurrences for the derivative quantities.

\subsection{Derivatives of modal Green's functions}
\label{sec:derivatives}

Following \cite{matviyenko1995azimuthal}, we introduce an alternative parametrization of (\ref{eq:mgf4}) in which the derivative structure becomes particularly simple. Define
\begin{align}
	a = R_0^2 \,, \quad  \mbox{and}\quad  b = \alpha R^2_0 = 2rr' \,,
\end{align}
so that (\ref{eq:mgf4}) becomes
\begin{align}
	G_m = \frac{1}{4\pi^2} \int_{0}^{\pi} \frac{e^{\I k  \sqrt{a-b\cos\theta}}}{\sqrt{a-b\cos\theta}} \cdot \cos(m\theta)\,\D \theta\,.
	\label{eq:mgf3}
\end{align}
By the chain rule, all derivatives of $G_m$ with respect to $r,r',z,z'$ can be written in terms of those with respect to $a$ and $b$\,. For example,
\begin{align}
	\frac{\partial G_m }{\partial r} &= 2r \frac{\partial G_m }{\partial a} + 2r'\frac{\partial G_m }{\partial b}  \,,
	\label{eq:dgdr_naive} \\
	\frac{\partial^2 G_m}{\partial r^2}&=2 \frac{\partial G_m}{\partial a}+4 r^2 \frac{\partial^2 G_m}{\partial a^2}+8 r r^{\prime} \frac{\partial^2 G_m}{\partial a \partial b}+4 r^{\prime 2} \frac{\partial^2 G_m}{\partial b^2} \,,
	\label{eq:d2gdr2_naive}
\end{align}
with analogous expressions for the remaining derivatives (see Appendix~\ref{sec:appendix_derivatives}).
Thus, to compute all first- and second-order derivatives of $G_m$\,, it suffices to compute
\begin{align}
	\frac{\partial G_m}{\partial a}\,, \quad \frac{\partial G_m}{\partial b}\,, \quad \frac{\partial^2 G_m}{\partial a^2}\,,  \quad \frac{\partial^2 G_m}{\partial a\partial b} \quad \mbox{and}  \quad \frac{\partial^2 G_m}{\partial b^2}\,.
	\label{eq:ab_derivs_needed}
\end{align}

Differentiating (\ref{eq:mgf3}) under the integral sign gives the required $(a,b)$-derivative integrals. For convenience, we define
\begin{align}
	w(\theta) = a - b\cos\theta\,.
	\label{eq:w}
\end{align}
For derivatives with respect to $a$, differentiating under the integral sign gives
\begin{align}
	\frac{\partial G_m}{\partial a} &=
	 \frac{1}{8\pi^2} \int_{0}^{\pi} \frac{e^{\I k \sqrt{w}}}{w} \sb{\I k - \frac{1}{\sqrt{w}}} \cos(m\theta)\,\D \theta\,,
	\label{eq:dgda} \\
	\frac{\partial^2 G_m}{\partial a^2} &=
	 \frac{1}{16\pi^2} \int_{0}^{\pi} \frac{e^{\I k \sqrt{w}}}{w} \sb{-\frac{k^2}{\sqrt{w}} - \frac{3\I k}{w}+ \frac{3}{w^{3/2}}} \cos(m\theta)\,\D \theta\,.
	\label{eq:d2gda2}
\end{align}
The derivatives with respect to $b$ have the same integrands, differing only by a factor of $-\cos\theta$ arising from $\frac{\partial w}{\partial b} = -\cos\theta$\,:
\begin{align}
	\frac{\partial G_m}{\partial b} &=
	 -\frac{1}{8\pi^2} \int_{0}^{\pi} \frac{e^{\I k \sqrt{w}}}{w} \sb{\I k - \frac{1}{\sqrt{w}}} \cos\theta\cos(m\theta)\,\D \theta\,,
	\label{eq:dgdb} \\
	\frac{\partial^2 G_m}{\partial a\partial b} &=
	 -\frac{1}{16\pi^2} \int_{0}^{\pi} \frac{e^{\I k \sqrt{w}}}{w} \sb{-\frac{k^2}{\sqrt{w}} - \frac{3\I k}{w}+ \frac{3}{w^{3/2}}} \cos\theta\cos(m\theta)\,\D \theta\,.
	\label{eq:d2gdadb}
\end{align}

In the near-singular regime where $\alpha\to 1^-$\,, the integrands develop a strong peak near $\theta = 0$\,, where $\cos\theta\approx 1$\,. The integrals are dominated by the contribution near the peak.
Consequently, $\frac{\partial G_m}{\partial a}$ and $-\frac{\partial G_m}{\partial b}$ are nearly equal in magnitude. Computing (\ref{eq:dgdr_naive}) directly leads to catastrophic cancellation since $r \approx r'$. The same cancellation occurs between $\frac{\partial^2 G_m}{\partial a^2}$ and $-\frac{\partial^2 G_m}{\partial a\partial b}$ in (\ref{eq:d2gdr2_naive})\,.
To circumvent this, we compute the combined quantities $\frac{\partial G_m}{\partial a} + \frac{\partial G_m}{\partial b}$ and $\frac{\partial^2 G_m}{\partial a^2} + \frac{\partial^2 G_m}{\partial a\partial b}$ directly:
\begin{align}
	\frac{\partial G_m}{\partial a} + \frac{\partial G_m}{\partial b} &=
	 \frac{1}{8\pi^2} \int_{0}^{\pi} \frac{e^{\I k \sqrt{w}}}{w} \sb{\I k - \frac{1}{\sqrt{w}}} (1-\cos\theta)\cos(m\theta)\,\D \theta\,,
	\label{eq:dgdapb} \\
	\frac{\partial^2 G_m}{\partial a^2} + \frac{\partial^2 G_m}{\partial a\partial b} &= \frac{1}{16\pi^2} \int_{0}^{\pi} \frac{e^{\I k \sqrt{w}}}{w} \sb{-\frac{k^2}{\sqrt{w}} - \frac{3\I k}{w}+ \frac{3}{w^{3/2}}} (1-\cos\theta)\cos(m\theta)\,\D \theta\,.
	\label{eq:d2gda2pab}
\end{align}
The factor $(1-\cos\theta)$ vanishes at $\theta = 0$\,, naturally removing the leading near-singular behavior of the integrand. The resulting integrals are orders of magnitude smaller as $\alpha \to 1^-$\,, and catastrophic cancellation is avoided.

With the combined quantities in hand, we define $\Delta r = r - r'$ and regroup the formulas (\ref{eq:dgdr_naive}) and (\ref{eq:d2gdr2_naive}) as
\begin{align}
	\frac{\partial G_m}{\partial r} &= 2\Delta r\, \frac{\partial G_m}{\partial a} + 2r' \bb{\frac{\partial G_m}{\partial a} + \frac{\partial G_m}{\partial b}}\,,
	\label{eq:dgdr_stable} \\
	\frac{\partial^2 G_m}{\partial r^2} &= 4r^2 \Big(\frac{\partial^2 G_m}{\partial a^2} + \frac{\partial^2 G_m}{\partial a\partial b}\Big) + 4r'^2 \Big(\frac{\partial^2 G_m}{\partial a\partial b} + \frac{\partial^2 G_m}{\partial b^2}\Big) - 4(\Delta r)^2 \frac{\partial^2 G_m}{\partial a\partial b} + 2\frac{\partial G_m}{\partial a}\,.
	\label{eq:d2gdr2_stable}
\end{align}
In (\ref{eq:dgdr_stable}), the potentially large quantities $\frac{\partial G_m}{\partial a}$ and $\frac{\partial G_m}{\partial b}$ appear only through their sum, while the remaining term is suppressed by the factor $\Delta r$\,. Similarly, in (\ref{eq:d2gdr2_stable}), the sums in the parentheses are also computed directly, and the only individually large term $\frac{\partial^2 G_m}{\partial a\partial b}$ is multiplied by $(\Delta r)^2$\,. The analogous numerically stable forms for all other cylindrical coordinate derivatives are collected in Appendix~\ref{sec:appendix_derivatives}.

Although (\ref{eq:ab_derivs_needed}) lists the five $(a,b)$-derivatives needed by the chain rule, in the stable form (\ref{eq:dgdr_stable})--(\ref{eq:d2gdr2_stable}) and the analogous formulas in Appendix~\ref{sec:appendix_derivatives}, $\frac{\partial^2 G_m}{\partial b^2}$ enters only through the sum $\frac{\partial^2 G_m}{\partial a\partial b}+\frac{\partial^2 G_m}{\partial b^2}$\,. The stable form therefore requires the seven quantities
\begin{align}
	&\frac{\partial G_m}{\partial a}\,, \quad \frac{\partial G_m}{\partial a} + \frac{\partial G_m}{\partial b}\,, \quad \frac{\partial^2 G_m}{\partial a^2}\,, \quad \frac{\partial^2 G_m}{\partial a^2} + \frac{\partial^2 G_m}{\partial a\partial b}\,, \label{eq:directquants}\\
	\intertext{and}
	&\frac{\partial G_m}{\partial b}\,, \quad \frac{\partial^2 G_m}{\partial a\partial b}\,, \quad \frac{\partial^2 G_m}{\partial a\partial b} + \frac{\partial^2 G_m}{\partial b^2}\,. \label{eq:recquants}
\end{align}
In Section~\ref{sec:recurrence}, we will show that the three quantities in (\ref{eq:recquants}) can be recovered via recurrence relations across Fourier modes, provided the four quantities in (\ref{eq:directquants}) have already been computed. As a result, it suffices in practice to evaluate only the four quantities in (\ref{eq:directquants}) directly --- we refer to these as the direct quantities --- while the three remaining quantities in (\ref{eq:recquants}) follow from the recurrences, and all derivatives of $G_m$ with respect to the cylindrical coordinates $(r,z,r',z')$ are then assembled via the chain rule of Appendix~\ref{sec:appendix_derivatives}.

\subsection{Recurrence relations}
\label{sec:recurrence}

In this section, we introduce recurrence relations satisfied by the modal Green's functions $G_m$ and their derivatives with respect to $a$ and $b$\,. The five-term recurrence for $G_m$ and the two-term recurrence for $\frac{\partial G_m}{\partial a}$ are due to \cite{matviyenko1995azimuthal}; the remaining relations for the second-order and combined derivatives are obtained by straightforward differentiation. The relations split into four recurrences in $m$, used to compute the four direct quantities in (\ref{eq:directquants}), and three pointwise identities that recover the remaining three quantities in (\ref{eq:recquants}).

\subsubsection{Five-term recurrence for $G_m$ \label{sec:recurrence5}}

The modal Green's functions satisfy the five-term recurrence relation \cite{matviyenko1995azimuthal}
\begin{align}
	c_0^m\, G_m + c_1^m\, G_{m+1} + c_{-1}^m\, G_{m-1} + c_2^m\, G_{m+2} + c_{-2}^m\, G_{m-2} = 0\,,
	\label{eq:fiveterm}
\end{align}
for $m \ge 2$\,, where the coefficients are given by
\begin{align}
	c_0^m &= 1 - \frac{k^2 b^2}{8a(m^2-1)} = 1 - \frac{\alpha^2\kappa^2}{8(m^2-1)}\,, \notag\\
	c_{\pm 1}^m &= -\frac{b}{a}\cdot\frac{2m\pm 1}{4m} = -\alpha\frac{2m\pm 1}{4m}\,, \notag\\
	c_{\pm 2}^m &= \frac{k^2 b^2}{16am(m\pm 1)} = \frac{\alpha^2\kappa^2}{16m(m\pm 1)}\,.
	\label{eq:fiveterm_coeffs}
\end{align}
Differentiating (\ref{eq:fiveterm}) with respect to $a$ and simplifying the right-hand side using (\ref{eq:fiveterm}) itself yields
\begin{align}
	c_0^m\, \frac{\partial G_m}{\partial a} + c_1^m\, \frac{\partial G_{m+1}}{\partial a} + c_{-1}^m\, \frac{\partial G_{m-1}}{\partial a} + c_2^m\, \frac{\partial G_{m+2}}{\partial a} + c_{-2}^m\, \frac{\partial G_{m-2}}{\partial a} = -\frac{G_m}{a}\,,
	\label{eq:fiveterm_da}
\end{align}
with the same coefficients $c_j^m$ as in (\ref{eq:fiveterm}). The two recurrences thus share an identical five-term structure and differ only in their right-hand sides.

The five-term recurrence (\ref{eq:fiveterm}) cannot be used directly as a forward or backward recursion to obtain all $G_m$\,, since there exist regimes in $m$ where both directions are unstable \cite{matviyenko1995azimuthal, garritano2022efficient}. A numerically stable use of (\ref{eq:fiveterm}) is deferred to Section~\ref{sec:allmodes}.

\subsubsection{Recurrences for derivatives \label{sec:recurrenced}}

For $b > 0$\,, the derivatives $\frac{\partial G_m}{\partial a}$ satisfy the two-term recurrence\,\cite{matviyenko1995azimuthal}
\begin{align}
	\frac{\partial G_{m+1}}{\partial a} - \frac{\partial G_{m-1}}{\partial a} = \frac{2m}{b}\, G_m\,, \quad m \ge 1\,.
	\label{eq:dgda_recurrence}
\end{align}
Since the integrands of $\frac{\partial G_m}{\partial b}$ and $\frac{\partial G_m}{\partial a}$ differ by a factor of $-\cos\theta$ (see Section~\ref{sec:derivatives}), and $\cos\theta \cdot e^{-\I m\theta} = \frac{1}{2}\bb{e^{-\I(m-1)\theta} + e^{-\I(m+1)\theta}}$\,, it follows that
\begin{align}
	\frac{\partial G_m}{\partial b} = -\frac{1}{2}\bb{\frac{\partial G_{m+1}}{\partial a} + \frac{\partial G_{m-1}}{\partial a}}\,, \quad m \ge 1\,.
	\label{eq:dgdb_recurrence}
\end{align}

For $b>0$\,, differentiating (\ref{eq:dgda_recurrence}) with respect to $a$ and applying the argument that gave (\ref{eq:dgdb_recurrence}) yields the analogous relations
\begin{align}
	\frac{\partial^2 G_{m+1}}{\partial a^2} - \frac{\partial^2 G_{m-1}}{\partial a^2} = \frac{2m}{b}\, \frac{\partial G_m}{\partial a}\,, \quad m \ge 1\,,
	\label{eq:d2gda2_recurrence}
\end{align}
and 
\begin{align}
	\frac{\partial^2 G_m}{\partial a\partial b} = -\frac{1}{2}\bb{\frac{\partial^2 G_{m+1}}{\partial a^2} + \frac{\partial^2 G_{m-1}}{\partial a^2}}\,, \quad m \ge 1\,.
	\label{eq:d2gdadb_recurrence}
\end{align}

Combining (\ref{eq:dgda_recurrence}) and (\ref{eq:dgdb_recurrence}) gives a recurrence for the cancellation-free first-derivative combination, valid for $b > 0$ and $m \ge 1$\,:
\begin{align}
	\bb{\frac{\partial G_{m+1}}{\partial a} + \frac{\partial G_{m+1}}{\partial b}} &= \bb{\frac{\partial G_{m-1}}{\partial a} + \frac{\partial G_{m-1}}{\partial b}} \notag\\
	&\quad + \frac{1}{b}\sb{-(m+1)\, G_{m+1} + 2m\, G_m - (m-1)\, G_{m-1}}\,.
	\label{eq:dgdapb_recurrence}
\end{align}
Differentiating this recurrence with respect to $a$ gives the corresponding recurrence for the second-order combined quantity:
\begin{align}
	\bb{\frac{\partial^2 G_{m+1}}{\partial a^2} + \frac{\partial^2 G_{m+1}}{\partial a\partial b}} &= \bb{\frac{\partial^2 G_{m-1}}{\partial a^2} + \frac{\partial^2 G_{m-1}}{\partial a\partial b}} \notag\\
	&\quad + \frac{1}{b}\sb{-(m+1)\, \frac{\partial G_{m+1}}{\partial a} + 2m\frac{\partial G_m}{\partial a} - (m-1)\frac{\partial G_{m-1}}{\partial a}}\,.
	\label{eq:d2gda2pab_recurrence}
\end{align}
Finally, applying the reasoning that gave (\ref{eq:dgdb_recurrence}) to $\frac{\partial^2 G_m}{\partial a^2} + \frac{\partial^2 G_m}{\partial a\partial b}$, we obtain
\begin{align}
	\frac{\partial^2 G_m}{\partial a\partial b} + \frac{\partial^2 G_m}{\partial b^2} = -\frac{1}{2}\Big[\Big(\frac{\partial^2 G_{m+1}}{\partial a^2} + \frac{\partial^2 G_{m+1}}{\partial a\partial b}\Big) + \Big(\frac{\partial^2 G_{m-1}}{\partial a^2} + \frac{\partial^2 G_{m-1}}{\partial a\partial b}\Big)\Big]\,,
	\label{eq:d2gdbpbb_recurrence}
\end{align}
for $m \ge 1$\,.

We observe that the four recurrences (\ref{eq:dgda_recurrence}), (\ref{eq:d2gda2_recurrence}), (\ref{eq:dgdapb_recurrence}), and (\ref{eq:d2gda2pab_recurrence}) do not amplify rounding errors, since each of them has the form $x_{m+1} = x_{m-1} + f_m$ with a known right-hand side $f_m$\,, whose solution is
\begin{align}
	x_{2k} = x_0 + \sum_{j=1}^{k} f_{2j-1}\,, \qquad x_{2k+1} = x_1 + \sum_{j=1}^{k} f_{2j}\,, \qquad k \ge 1\,,
	\label{eq:twoterm_solution}
\end{align}
simply the starting value plus a sum of the right-hand sides of the same parity; recursing downward from starting values $x_M$\,, $x_{M-1}$ gives likewise
\begin{align}
	x_{M-2k} = x_M - \sum_{j=1}^{k} f_{M-2j+1}\,, \qquad x_{M-1-2k} = x_{M-1} - \sum_{j=1}^{k} f_{M-2j}\,, \qquad k \ge 1\,.
	\label{eq:twoterm_solution_down}
\end{align}
In either direction a perturbation of any term enters with coefficient one, so rounding errors accumulate additively, at most linearly in the number of steps, and are never magnified by a growing homogeneous solution. Consequently, they can be applied as stable upward recurrences from starting values at $m = 0, 1$\,. In the decay regime, the upward recurrence loses relative accuracy beyond the transition point $m^*$ (see Section~\ref{sec:regimes}): for $m > m^*$\,, the partial sums in (\ref{eq:twoterm_solution}) are of size $O(1)$ while $x_m$ is exponentially small, so each addition commits an absolute rounding error of the order of machine epsilon, far larger than the result. In (\ref{eq:twoterm_solution_down}), by contrast, the terms are added from the smallest to the largest: the partial sums remain comparable in magnitude to the most recently added terms, so the rounding error of each addition is relative to the partial sum, and $x_m$ is obtained to relative accuracy. Therefore, a downward recurrence starting from some $m > m^*$ should be used instead. The three relations (\ref{eq:dgdb_recurrence}), (\ref{eq:d2gdadb_recurrence}), and (\ref{eq:d2gdbpbb_recurrence}) are local identities and are always stable.

Together with starting values at $m = 0, 1$\,, these relations determine all derivatives once $G_m$\,, $m = 0, 1, \ldots, M$\,, is known: the four direct quantities in (\ref{eq:directquants}) are obtained by the recurrences (\ref{eq:dgda_recurrence}), (\ref{eq:d2gda2_recurrence}), (\ref{eq:dgdapb_recurrence}), (\ref{eq:d2gda2pab_recurrence}); the three remaining quantities in (\ref{eq:recquants}) follow from the pointwise identities (\ref{eq:dgdb_recurrence}), (\ref{eq:d2gdadb_recurrence}), (\ref{eq:d2gdbpbb_recurrence}).

\begin{remark}
	The $b>0$ restriction in (\ref{eq:dgda_recurrence}), (\ref{eq:d2gda2_recurrence}), (\ref{eq:dgdapb_recurrence}), and (\ref{eq:d2gda2pab_recurrence}) does not impose any limitation. When $b = 2rr' \approx 0$\,, we have $\alpha = b/a \approx 0$\,, which is precisely the near-axis regime (Section~\ref{sec:regimes}) where the power series expansion (\ref{eq:nearaxis}) and its derivatives apply.
\end{remark}

\section{Contour deformation for a single Fourier mode}
\label{sec:contour}

In this section we evaluate a single modal Green's function $G_m$ via contour deformation, building on the approach of \cite{garritano2022efficient} and extending it to the derivative quantities of Section~\ref{sec:derivatives}. The result is a uniformly accurate, $\kappa$-independent quadrature for $G_m$ at any single mode $m$. The same deformed contour applies to all first- and second-order derivative integrals. For the near-singular integral on $\gamma_1$, we replace the monomial recurrence used in \cite{gustafsson2010accurate,garritano2022efficient} with precomputed Generalized Gaussian Quadratures \cite{bremer2010nonlinear} tailored to the singular factors arising in the derivative integrands. As in \cite{garritano2022efficient}, the cost of evaluating a single $G_m$ is $O(m)$, independent of the wavenumber and source-target separation.

\subsection{Integrals to be evaluated}
\label{sec:integrals_eval}
We first rewrite the four direct quantities in (\ref{eq:directquants}) in the $(\alpha, \kappa)$ variables used by the contour deformation. To match the notation of \cite{garritano2022efficient}, we convert the integral representations from the $(a,b)$-parametrization of Section~\ref{sec:derivatives} to the one introduced in Section~\ref{sec:parametrization}.

Since $w = a - b\cos\theta = R_0^2(1 - \alpha\cos\theta)$ from (\ref{eq:w})\,, we have $\sqrt{w} = R_0\sqrt{1-\alpha\cos\theta}$ and $k\sqrt{w} = \kappa\sqrt{1-\alpha\cos\theta}$\,. The integrals of Section~\ref{sec:derivatives} can therefore be rewritten entirely in terms of $\alpha$ and $\kappa$\,. We apply the substitution $x = \cos\theta$ and define
\begin{align}
	s(x) = \sqrt{1-\alpha x}\,.
	\label{eq:sx}
\end{align}
 Since the Chebyshev polynomial of the first kind satisfies $\cos(m\theta) = T_m(x)$\,, the modal Green's function in (\ref{eq:mgf4}) becomes
\begin{align}
	G_m = \frac{1}{4\pi^2 R_0} \int_{-1}^{1} \frac{e^{\I\kappa s}}{s} \cdot \frac{T_m(x)}{\sqrt{1-x^2}}\,\D x\,.
	\label{eq:gm_cheb}
\end{align}
Similarly, the derivatives in (\ref{eq:dgda})--(\ref{eq:d2gda2}) become
\begin{align}
	\frac{\partial G_m}{\partial a} &= \frac{1}{8\pi^2 R_0^3} \int_{-1}^{1} e^{\I\kappa s}\bb{\frac{\I\kappa}{s^2} - \frac{1}{s^3}} \frac{T_m(x)}{\sqrt{1-x^2}}\,\D x\,,
	\label{eq:dgda_cheb} \\
	\frac{\partial^2 G_m}{\partial a^2} &= \frac{1}{16\pi^2 R_0^5} \int_{-1}^{1} e^{\I\kappa s}\bb{-\frac{\kappa^2}{s^3} - \frac{3\I\kappa}{s^4} + \frac{3}{s^5}} \frac{T_m(x)}{\sqrt{1-x^2}}\,\D x\,,
	\label{eq:d2gda2_cheb}
\end{align}
and the combined quantities in (\ref{eq:dgdapb})--(\ref{eq:d2gda2pab}) become
\begin{align}
	\frac{\partial G_m}{\partial a} + \frac{\partial G_m}{\partial b} &= \frac{1}{8\pi^2 R_0^3} \int_{-1}^{1} e^{\I\kappa s}\bb{\frac{\I\kappa}{s^2} - \frac{1}{s^3}} \frac{(1-x)\,T_m(x)}{\sqrt{1-x^2}}\,\D x\,,
	\label{eq:dgdapb_cheb} \\
	\frac{\partial^2 G_m}{\partial a^2} + \frac{\partial^2 G_m}{\partial a\partial b} &= \frac{1}{16\pi^2 R_0^5} \int_{-1}^{1} e^{\I\kappa s}\bb{-\frac{\kappa^2}{s^3} - \frac{3\I\kappa}{s^4} + \frac{3}{s^5}} \frac{(1-x)\,T_m(x)}{\sqrt{1-x^2}}\,\D x\,.
	\label{eq:d2gda2pab_cheb}
\end{align}

\subsection{Steepest descent contours for real $k$}
\label{sec:sdcontours}

In the following, we describe the contour deformation used to evaluate the integrals (\ref{eq:gm_cheb})--(\ref{eq:d2gda2pab_cheb}) for real $k \ge 0$\,, following the approach of \cite{garritano2022efficient}. On these contours, the integrals can be computed accurately with $O(m)$ quadrature points, independent of $\kappa$\,. We state only the key results here and refer the reader to \cite{garritano2022efficient} for the analysis and implementation details.

The integrands in (\ref{eq:gm_cheb})--(\ref{eq:d2gda2pab_cheb}) have branch points at $z = \pm 1$ from $\sqrt{1-z^2}$ and at $z = 1/\alpha$ from $s(z) = \sqrt{1-\alpha z}$; the derivative integrands have stronger singularities at $z = 1/\alpha$. Away from these singularities and the associated branch cuts, the integrands are analytic. We take $\sqrt{1-z^2} = \sqrt{1-z}\,\sqrt{1+z}$\,, with both factors on their principal branches, giving cuts $(-\infty,-1]$ and $[1,\infty)$\,. This branch is positive on $(-1,1)$ and extends analytically into the lower half-plane. We likewise take $s(z)$ on its principal branch, with cut $[1/\alpha,\infty)$\,. The deformed contour lies in the lower half-plane except at its endpoints $z = \pm 1$ and crosses none of these cuts. By Cauchy's theorem, the integral over $[-1,1]$ can be deformed into the lower half of the complex plane. The steepest descent contours $\gamma_1$ and $\gamma_2$\,, emanating from $z = 1$ and $z = -1$ respectively, are chosen so that the spherical wave term $e^{\I\kappa s}$ is non-oscillatory and decays exponentially. Using the notation of (\ref{eq:gmparas}), the contours are parametrized by
\begin{align}
	\gamma_1(\tau) = 1 + \tau^4 - 2\I\beta_-\tau^2\,, \quad \gamma_2(\tau) = -1 + \tau^4 - 2\I\beta_+\tau^2\,, \quad \tau \ge 0\,.
	\label{eq:gamma12}
\end{align}
On $\gamma_1$\,, it is easy to verify that $s(\gamma_1(\tau)) = \I\sqrt{\alpha}\,(\tau^2 - \I\beta_-)$\,, so that
\begin{align}
	e^{\I\kappa s(\gamma_1(\tau))} = e^{-\kappa\sqrt{\alpha}\,(\tau^2 - \I\beta_-)}\,,
	\label{eq:decay1}
\end{align}
which decays as $e^{-\kappa\sqrt{\alpha}\,\tau^2}$ for $\tau > 0$\,. An analogous formula holds on $\gamma_2$ with $\beta_-$ replaced by $\beta_+$\,.

While the spherical wave term decays on $\gamma_1$ and $\gamma_2$\,, the Chebyshev polynomial $T_m(z)$ grows exponentially along these contours. To avoid catastrophic cancellation, the contours are truncated at their intersection with a Bernstein ellipse $E_\rho$ in the complex plane, parametrized by
\begin{align}
	E_\rho(\theta) = a_e\cos\theta - \I b_e\sin\theta\,, \quad \theta \in [0, 2\pi)\,,
	\label{eq:bernstein}
\end{align}
where $a_e = \frac{1}{2}(\rho + \rho^{-1})$\,, $b_e = \frac{1}{2}(\rho - \rho^{-1})$\,, and $\rho > 1$\,. We observe that (\ref{eq:bernstein}) parametrizes the ellipse clockwise; we adopt this convention because the steepest descent contours lie in the lower half-plane. On the ellipse, the size of the Chebyshev polynomial satisfies the following bound:
\begin{align}
	|T_m(z)| \le \frac{1}{2}(\rho^m + \rho^{-m}) \le \rho^m\,, \quad z \in E_\rho\,.
	\label{eq:tm_bound}
\end{align}
The parameter $\rho$ is therefore chosen so that the growth of $T_m$ on the ellipse remains bounded.

The contour $\gamma_1$ intersects $E_\rho$ at $\theta = \theta_1$ (near $z = 1$), and $\gamma_2$ intersects $E_\rho$ at $\theta = \theta_2$ (near $z = -1$). Matching real and imaginary parts of $\gamma_1(\tau_1) = E_\rho(\theta_1)$ and $\gamma_2(\tau_2) = E_\rho(\theta_2)$ leads to quadratic equations (see \cite{garritano2022efficient}), whose roots satisfy
\begin{align}
	\cos\theta_1 = \frac{-a_e + \sqrt{a_e^2 + \frac{b_e^2}{\beta_-^2}\bb{1 + \frac{b_e^2}{4\beta_-^2}}}}{\frac{b_e^2}{2\beta_-^2}}\,, \quad
	\cos\theta_2 = \frac{-a_e + \sqrt{a_e^2 - \frac{b_e^2}{\beta_+^2}\bb{1 - \frac{b_e^2}{4\beta_+^2}}}}{\frac{b_e^2}{2\beta_+^2}}\,.
	\label{eq:costheta12}
\end{align}
Then we solve for $\theta_1$ and $\theta_2$, from which the contour intersections in the $\tau$ parametrization are obtained by the formulas
\begin{align}
	\tau_1 = \sqrt{\frac{b_e\sin\theta_1}{2\beta_-}}\,, \quad \tau_2 = \sqrt{\frac{b_e\sin\theta_2}{2\beta_+}}\,.
	\label{eq:tau12}
\end{align}

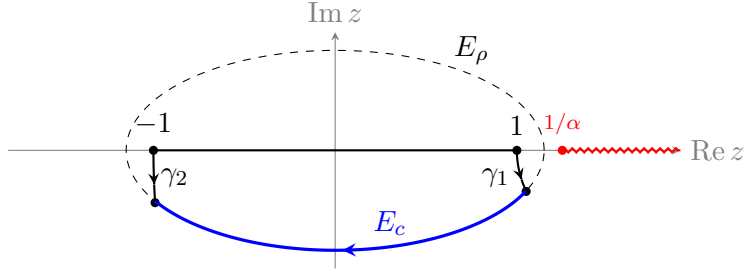
\begin{figure}[ht]
\centering
\begin{tikzpicture}[scale=2.4, >=stealth]
% axes
\draw[gray, ->] (-1.8,0) -- (1.9,0) node[right] {$\mathrm{Re}\, z$};
\draw[gray, ->] (0,-0.65) -- (0,0.65) node[above] {$\mathrm{Im}\, z$};

% real interval [-1,1]
\draw[thick] (-1,0) -- (1,0);
\fill (-1,0) circle (0.025) node[above=2pt] {$-1$};
\fill (1,0) circle (0.025) node[above=2pt] {$1$};

% branch point
\fill[red] (1.25,0) circle (0.025) node[above=2pt] {$\scriptstyle 1/\alpha$};
% branch cut
\draw[red, thick, decorate, decoration={zigzag, segment length=3pt, amplitude=0.8pt}] (1.25,0) -- (1.9,0);

% Bernstein ellipse (full, dashed)
\draw[dashed] plot[domain=0:360, samples=200, smooth]
  ({1.15*cos(\x)},{0.55*sin(\x)});
\node at (0.75,0.55) {$E_\rho$};

% gamma_1: from (1,0) downward to intersection (~1.051, -0.226)
% Im = -t^2 (reflected from upper half-plane)
\draw[thick, ->] plot[domain=0:0.40, samples=50, smooth, variable=\t]
  ({1+\t*\t*\t*\t},{-\t*\t});
\draw[thick] plot[domain=0.40:0.475, samples=20, smooth, variable=\t]
  ({1+\t*\t*\t*\t},{-\t*\t});
\node[right] at (0.75,-0.15) {$\gamma_1$};

% gamma_2: from (-1,0) downward to intersection (~-0.991, -0.288)
\draw[thick, ->] plot[domain=0:0.25, samples=50, smooth, variable=\t]
  ({-1+\t*\t*\t*\t},{-3*\t*\t});
\draw[thick] plot[domain=0.25:0.31, samples=20, smooth, variable=\t]
  ({-1+\t*\t*\t*\t},{-3*\t*\t});
\node[left] at (-0.75,-0.15) {$\gamma_2$};

% intersection points
\fill (1.051, -0.226) circle (0.025);
\fill (-0.991, -0.288) circle (0.025);

% E_c arc (lower portion of ellipse between intersection points)
% -theta_1 = -24 deg = 336 deg, -theta_2 = -149.5 deg = 210.5 deg
\draw[blue, very thick] plot[domain=210.5:336, samples=100, smooth]
  ({1.15*cos(\x)},{0.55*sin(\x)});
\node[blue] at (0.3,-0.40) {$E_c$};

% direction arrow on E_c (going from near 1 toward near -1, i.e., theta_1 to theta_2)
\draw[blue, very thick, ->] plot[domain=274:272, samples=2, smooth]
  ({1.15*cos(\x)},{0.55*sin(\x)});

\end{tikzpicture}
\caption{Contour geometry in the complex $z$-plane for $\alpha = 0.8$ ($\beta_- = 0.5$, $\beta_+ = 1.5$). The integration path along $[-1,1]$ is deformed into the steepest descent contours $\gamma_1$\,, $\gamma_2$ and the Bernstein ellipse arc $E_c$ in the lower half-plane. The branch point at $z = 1/\alpha$ and the branch cut (zigzag) are shown in red.}
\label{fig:contours}
\end{figure}

After defining the contour, the integral is decomposed as
\begin{align}
	\int_{-1}^{1} = -\int_{\gamma_1} + \int_{\gamma_2} - \int_{E_c}\,,
	\label{eq:contour_decomp}
\end{align}
where $\gamma_1$ and $\gamma_2$ are both traversed from $z = \pm 1$ downward to the ellipse (i.e., $\tau = 0$ to $\tau_1$ and $\tau_2$ respectively), and $E_c$ denotes the arc of the Bernstein ellipse from $\theta_1$ to $\theta_2$\,. The contour geometry is illustrated in Figure~\ref{fig:contours}.

\begin{remark}
	The ellipse parameter $\rho$ is chosen as $\rho = B^{1/m}$ for a prescribed constant $B$\,, so that $|T_m(z)| \le \rho^m = B$ on $E_\rho$ by (\ref{eq:tm_bound}). This means that the cancellation between the growing Chebyshev polynomial and the decaying spherical wave on $\gamma_1$ and $\gamma_2$ costs at most $\log_{10} B$ digits of accuracy. In this paper, we choose $B=100$ so that we lose at most two digits due to cancellation. The choice of $B$ is a trade-off: $\rho$ controls both the size of $T_m$ on the ellipse and the size of the ellipse itself. A larger $B$ costs more digits to cancellation, but the larger ellipse keeps the arc $E_c$ farther from the singularities of the integrand at $z = \pm 1$ and $z = 1/\alpha$\,, so fewer quadrature nodes are needed there; a smaller $B$ does the opposite. Quantitatively, it is shown in \cite{garritano2022efficient} that attaining the accuracy $B\epsilon$ on $E_c$\,, with $\epsilon$ the machine precision, requires approximately $m \log(2/\epsilon) / (2\log B)$ Gauss--Legendre nodes. For $B = 100$\,, this is about $4m$ nodes in double precision, consistent with the conservative choice of $5m$ used in Section~\ref{sec:quadrature}.
	\label{rmk:rhob}
\end{remark}

\begin{remark}
	\label{rem:lowmode}
	For $m = 0$\,, $T_0(z) = 1$ does not grow in the complex plane, so the Bernstein ellipse truncation is unnecessary and the steepest descent contours $\gamma_1$\,, $\gamma_2$ extend indefinitely into the lower half-plane; this is the setting of \cite{gustafsson2010accurate}. For small $m > 0$\,, the ellipse parameter $\rho = B^{1/m}$ is large and the contours are truncated only far from the real axis. For a more uniform treatment, we use the contour geometry of $m = 5$ for all modes $m \le 4$\,: the Bernstein ellipse, intersection points, and truncation parameters are computed with $m = 5$\,, while the integrand is evaluated with the correct mode number $m$\,.
\end{remark}

\subsection{Integrands on each contour segment}
\label{sec:contour_integrands}

This section contains the integrand formulas on each contour segment. We introduce the notation $u_- = \tau^2 - \I\beta_-$ and $u_+ = \tau^2 - \I\beta_+$\,, so that $\sqrt{1-\alpha\gamma_1(\tau)} = \I\sqrt{\alpha}\, u_-$ and $\sqrt{1-\alpha\gamma_2(\tau)} = \I\sqrt{\alpha}\, u_+$ (see (\ref{eq:decay1})). We apply the decomposition (\ref{eq:contour_decomp}) to the modal Green's function (\ref{eq:gm_cheb}) to obtain
\begin{align}
	G_m = \frac{1}{4\pi^2 R_0}\bb{-I_{\gamma_1}^{(0)} + I_{\gamma_2}^{(0)} - I_{E_c}^{(0)}}\,, 
	\label{eq:gm_decomp}
\end{align}
where the three integrals on the contours $\gamma_1$\,, $\gamma_2$\,, and $E_c$ are respectively given by the formulas
\begin{align}
	I_{\gamma_1}^{(0)} = -\frac{4}{\sqrt{\alpha}} \int_{0}^{\tau_1} \frac{e^{-\kappa\sqrt{\alpha}\, u_-}}{\sqrt{\tau^2 - 2\I\beta_-}\,\sqrt{1+\gamma_1(\tau)}}\, T_m\bb{\gamma_1(\tau)}\, \D\tau\,, 
	\label{eq:igamma1_gm}
\end{align}
\begin{align}
	I_{\gamma_2}^{(0)} = -\frac{4\I}{\sqrt{\alpha}} \int_{0}^{\tau_2} \frac{e^{-\kappa\sqrt{\alpha}\, u_+}}{\sqrt{\tau^2 - 2\I\beta_+}\,\sqrt{1-\gamma_2(\tau)}}\, T_m\bb{\gamma_2(\tau)}\, \D\tau\,,
	\label{eq:igamma2_gm}
\end{align}
and 
\begin{align}
	I_{E_c}^{(0)} = \int_{\theta_1}^{\theta_2} \frac{e^{\I\kappa\sqrt{1-\alpha z}}}{\sqrt{1-\alpha z}\,\sqrt{1-z^2}}\, T_m(z)\, z'(\theta)\,\D\theta\,,
	\label{eq:iec_gm}
\end{align}
with $z'(\theta) = -a_e\sin\theta - \I b_e\cos\theta$ and $T_m(z(\theta)) = \frac{1}{2}(\rho^m e^{-\I m\theta} + \rho^{-m} e^{\I m\theta})$ on the ellipse.

The derivative integrands from Section~\ref{sec:integrals_eval} are evaluated on the same three contour segments, with the factor $1/\sqrt{1-\alpha z}$ in the $G_m$ integrand replaced by the appropriate derivative factor. The explicit decompositions for all four derivative quantities $\frac{\partial G_m}{\partial a}$\,, $\frac{\partial G_m}{\partial a} + \frac{\partial G_m}{\partial b}$\,, $\frac{\partial^2 G_m}{\partial a^2}$\,, and $\frac{\partial^2 G_m}{\partial a^2} + \frac{\partial^2 G_m}{\partial a\partial b}$ are summarized in Appendix~\ref{sec:appendix_contour}.

\subsection{Quadrature on each contour segment}
\label{sec:quadrature}
In this section, we summarize the quadratures for accurately approximating (\ref{eq:igamma1_gm})--(\ref{eq:iec_gm}). Due to the contour choice of $\gamma_1$ and $\gamma_2$\,, the only oscillatory term in (\ref{eq:igamma1_gm}) and (\ref{eq:igamma2_gm}) is $T_m(z)$\,. As explained in Remark\,\ref{rmk:rhob}, we choose the ellipse parameter as $\rho = B^{1/m}$ with $B=100$\,. With this choice, it is proved in \cite{garritano2022efficient} that the ellipse geometry ensures that the integrand in (\ref{eq:igamma1_gm}) and (\ref{eq:igamma2_gm}) has a bounded number of oscillations. Consequently, a bounded number of quadrature nodes suffices on $\gamma_1$ and $\gamma_2$\,, independent of $m$\,.

%As explained in Remark\,\ref{rmk:rhob}, we choose the ellipse parameter as $\rho = B^{1/m}$ with $B=100$\,. With this choice, the Bernstein ellipse shrinks as $m$ increases since $b_e = \frac{1}{2}(\rho - \rho^{-1}) = O(1/m)$\,. The steepest descent contours $\gamma_1$ and $\gamma_2$ are truncated at their intersection with this ellipse, resulting in contour segments of arc length $O(1/m^2)$ in the complex plane. Since the oscillation density of $T_m$ near $z = \pm 1$ is $O(m^2)$\,, the total number of oscillations of $T_m$ on each steepest descent segment is $O(m^2 \cdot 1/m^2) = O(1)$\,. Consequently, a bounded number of quadrature nodes suffices on $\gamma_1$ and $\gamma_2$\,, independent of $m$\,. 

\subsubsection{Contour $\gamma_2$ (near $z = -1$)}

Since $\beta_+ \ge 1$\,, the integrand (\ref{eq:igamma2_gm}) on $\gamma_2$ is smooth on $[0, \tau_2]$ and oscillates $O(1)$ times, so standard Gauss--Legendre quadrature with a fixed number of nodes suffices. In this paper, we use 32 nodes for all $\gamma_2$ integrals, including those for the derivative quantities.

\subsubsection{Ellipse arc $E_c$}

On the ellipse arc $E_c$\,, the Chebyshev polynomial $T_m$ in (\ref{eq:iec_gm}) oscillates $O(m)$ times (since $\theta$ traverses an interval of length $O(1)$ and $T_m(E_\rho(\theta))$ involves $e^{\pm \I m\theta}$). Furthermore, the ellipse contour stays at distance $O(1/m^2)$ from $z = \pm 1$\,, where the integrand of (\ref{eq:iec_gm}) is singular.
As proved in \cite{garritano2022efficient}, this segment therefore requires $O(m)$ quadrature nodes. In our implementation, we choose $5m$ Gauss--Legendre nodes on $[\theta_1, \theta_2]$. The same argument applies to the derivative integrands, so the same number of nodes is used. This segment is the origin of the dominant $O(m)$ cost of evaluating a single mode for large $m$.
\subsubsection{Contour $\gamma_1$ (near $z = +1$)}

The integrand in (\ref{eq:igamma1_gm}) on $\gamma_1$ likewise oscillates $O(1)$ times. However, due to the factor $\frac{1}{\sqrt{\tau^2 - 2\I\beta_-}}$\,, when $\alpha \approx 1$ (so $\beta_-$ is small), the integrand is near-singular in the neighborhood of $\tau = 0$\,, preventing the use of smooth quadrature as on $\gamma_2$\,. In \cite{garritano2022efficient}, this is handled by expanding the smooth part of the integrand in monomials and evaluating the resulting singular integrals via a recurrence relation. For the derivative integrands, this monomial recurrence becomes more involved due to the additional factors of $u_-^{-1},\ldots,u_-^{-4}$ with $u_- = \tau^2 - \I\beta_-$ in the derivative integrands of Appendix~\ref{sec:appendix_contour}.

We instead use precomputed Generalized Gaussian Quadratures \cite{bremer2010nonlinear}. These quadratures integrate the singular factors directly, avoiding separate monomial recurrences for each derivative integrand. Every singular factor appearing in the $\gamma_1$ integrands of Section~\ref{sec:contour} and Appendix~\ref{sec:appendix_contour} can be represented as a polynomial multiple of one of the five families
\begin{equation}
\frac{p(\tau)}{(\tau^2 - \I\beta_-)^j\,\sqrt{\tau^2 - 2\I\beta_-}}\,, \qquad j = 0, 1, 2, 3, 4\,,
\label{eq:ggq_weights_core}
\end{equation}
for any polynomial $p(\tau)$ of degree at most $63$\,. We construct a Generalized Gaussian Quadrature on $[0,1]$ that integrates these five families exactly, with $\beta_-$ ranging over a dyadic decomposition of $[2^{-90}, 2^{-1}]$\,.

These five families are not independent of each other: each less singular family $j < 4$ lies in the polynomial span of the most singular family $j = 4$\,, by the trivial identity
\begin{equation}
\frac{p(\tau)}{(\tau^2 - \I\beta_-)^j\,\sqrt{\tau^2 - 2\I\beta_-}} = \frac{p(\tau)\,(\tau^2 - \I\beta_-)^{4-j}}{(\tau^2 - \I\beta_-)^4\,\sqrt{\tau^2 - 2\I\beta_-}}\,,
\end{equation}
which raises the polynomial degree from $N$ to $N + 2(4-j)$\,. Thus the five-family construction is essentially equivalent to constructing a quadrature for the most singular family at a slightly higher polynomial order. In practice, this increases the number of nodes only mildly compared with a quadrature constructed for $G_m$ alone (which corresponds to the $j = 0$ family of (\ref{eq:ggq_weights_core})).

In the well-separated regime $\beta_- > 2^{-1}$\,, the integrands are smooth and a $32$-node Gauss--Legendre quadrature suffices. For smaller $\beta_-$\,, our implementation selects one of $18$ precomputed Generalized Gaussian Quadratures according to the value of $\beta_-$\,, with node counts ranging from $67$ for the least singular case to $130$ for the most singular.\footnote{Available at \href{https://github.com/han-wen-zhang/modal_helmholtz_quadrature}{\mbox{github.com/}\allowbreak\mbox{han-wen-zhang/}\allowbreak\mbox{modal\_helmholtz\_quadrature}}.}

\begin{remark}
\label{rem:shared_cost}
A key advantage of using the same quadrature nodes for all quantities (\ref{eq:gm_cheb})--(\ref{eq:d2gda2pab_cheb}) is that the expensive parts of the function evaluations (the exponential, Chebyshev polynomial, and contour parametrization) can be shared across all five integrands. As a result, evaluating all of (\ref{eq:gm_cheb})--(\ref{eq:d2gda2pab_cheb}) simultaneously adds only about $20\%$ to the cost of evaluating $G_m$ alone.
\end{remark}

\begin{remark}
The overall cost of evaluating a single mode $G_m$ and its derivatives ranges from $64 + 5m$ quadrature points for well-separated cases and $162 + 5m$ for the most singular cases, independent of the frequency parameter $\kappa$, and hence of the wavenumber $k$. Thus the number of quadrature nodes is bounded uniformly in both frequency and source-target separation. 
These numbers are conservative and no attempt has been made to optimize them.
\end{remark}

\begin{remark}
\label{rem:complexk}
	Although we restrict to real $k \ge 0$ in this paper, the contour deformation extends to complex $k$ (e.g., attenuating media with $\Im{k} > 0$). As shown in \cite{garritano2022efficient}, the steepest descent contours $\gamma_1$ and $\gamma_2$ rotate in the complex plane. Their intersections with the Bernstein ellipse are obtained by Newton iteration at little additional cost, and the distribution of quadrature nodes among $\gamma_1$\,, $\gamma_2$ and $E_c$ varies with $\arg\kappa$\,, while the total remains $O(m)$\,. The quadrature rules on $\gamma_1$\,, tabulated in $\beta_-$\,, then acquire $\arg\kappa$ as a second parameter and are constructed in the same manner. For nearly imaginary $\kappa$\,, the spherical-wave term $e^{\I\kappa s} \approx e^{-|\kappa| s}$ is non-oscillatory, as for the modified Helmholtz equation, and the steepest descent contour $\gamma_1$ approaches the singularity at $z = -1$\,; in this case, the integration is carried out along the steepest descent contour associated with $\Re\kappa$\,, which is well separated from the singularity; see Section~4.4.4 of \cite{garritano2022efficient}.
\end{remark}

\section{Evaluation of all Fourier modes}
\label{sec:allmodes}

Section~\ref{sec:contour} gives an $O(m)$ contour-deformation method for evaluating a single mode $G_m$ and its derivatives. We now combine this method with the five-term recurrence (\ref{eq:fiveterm}) to compute all modes $G_0, \ldots, G_M$ in $O(M)$ time, for any $M \ge 5$\,. The five-term recurrence (\ref{eq:fiveterm}) cannot be used directly as a forward or backward recursion to obtain all $G_m$, since there exist regimes in $m$ where both directions are unstable \cite{matviyenko1995azimuthal, garritano2022efficient}. We therefore use it as a boundary-value problem in the mode index $m$\,: the boundary values $G_0, G_1, G_{M-1}, G_M$ are supplied by contour deformation, and the interior values are obtained from a banded linear system. We demonstrate experimentally that this formulation gives \emph{componentwise} accuracy even when the pentadiagonal system has a large condition number (by componentwise accuracy we mean that each individual computed mode $G_m$ has small relative error, rather than merely contributing to a small normwise error). It should be observed that the analysis of this phenomenon remains incomplete and is currently under vigorous investigation.

In the non-decay regime, $G_{M-1}$ and $G_M$ are evaluated accurately by the contour deformation of Section~\ref{sec:contour}. Enforcing (\ref{eq:fiveterm}) for $m = 2, 3, \ldots, M-2$ is then equivalent to solving the following pentadiagonal linear system of size $(M-3) \times (M-3)$\,:
{\setlength{\mathindent}{30pt}%
\begin{align}
	\begin{pmatrix}
		c_0^2 & c_1^2 & c_2^2 \\[2pt]
		c_{-1}^3 & c_0^3 & c_1^3 & c_2^3 \\[2pt]
		c_{-2}^4 & c_{-1}^4 & c_0^4 & c_1^4 & c_2^4 \\[4pt]
		 & \ddots & \ddots & \ddots & \ddots & \ddots \\[4pt]
		 & & c_{-2}^{M\!-\!3} & c_{-1}^{M\!-\!3} & c_0^{M\!-\!3} & c_1^{M\!-\!3} \\[2pt]
		 & & & c_{-2}^{M\!-\!2} & c_{-1}^{M\!-\!2} & c_0^{M\!-\!2}
	\end{pmatrix}
	\begin{pmatrix} G_2 \\ G_3 \\ G_4 \\ \vdots \\ G_{M\!-\!3} \\ G_{M\!-\!2} \end{pmatrix}
	=
	-\begin{pmatrix} c_{-2}^2 G_0 + c_{-1}^2 G_1 \\ c_{-2}^3 G_1 \\ 0 \\ \vdots \\ 0 \\ c_2^{M\!-\!3} G_{M\!-\!1} \\ c_1^{M\!-\!2} G_{M\!-\!1} + c_2^{M\!-\!2} G_M \end{pmatrix}\,,
	\label{eq:pentadiag}
\end{align}}%
where the known boundary values $G_0$\,, $G_1$\,, $G_{M-1}$\,, $G_M$ appear on the right-hand side. This system is solved in $O(M)$ time by a banded solver. We choose a banded QR factorization with Givens rotations \cite{golub2013matrix} for its simplicity: it requires no pivoting and is unconditionally backward stable, at the price of increasing the bandwidth of the matrix by two. We observe experimentally that banded LU with partial pivoting produces essentially the same accuracy at a smaller constant factor; since the dominant cost is the contour evaluation of $G_{M-1}$ and $G_M$, the choice has negligible effect on total runtime. 

As explained in Section\,\ref{sec:regimes}, when the requested mode range enters the decay regime, evaluating $G_{M-1}$\,, $G_M$ via the integral representation loses relative accuracy. In this case, we instead choose some $M'$\,, where $\abs{G_{M'-1}}$ and $\abs{G_{M'}}$ are well below machine epsilon, set $G_{M'-1}$ and $G_{M'}$ to zero on the right-hand side of (\ref{eq:pentadiag}), and solve the $(M'-3) \times (M'-3)$ linear system. This is the five-term analogue of Miller's algorithm for computing the minimal solution of a recurrence. 

\begin{remark}
	The idea of casting a recurrence as a linear system appears in Gautschi's survey \cite{gautschi1967computational}, where it is dismissed as rarely practical because the solution values at high mode numbers are typically unavailable. In our setting, the contour deformation of Section~\ref{sec:contour} supplies those values directly. The approach was subsequently developed by Olver \cite{olver1967numerical}, who formulates the three-term recurrence as a tridiagonal boundary-value system for the minimal solution (as in our decay-regime approach above) and gives a rigorous relative-error analysis. The pentadiagonal boundary-value problem arising here appears to be less standard, and we are not aware of an analogous relative-error analysis.
\end{remark}

\begin{figure}[!ht]
\centering
\begin{subfigure}[b]{0.48\textwidth}
\includegraphics[width=\textwidth]{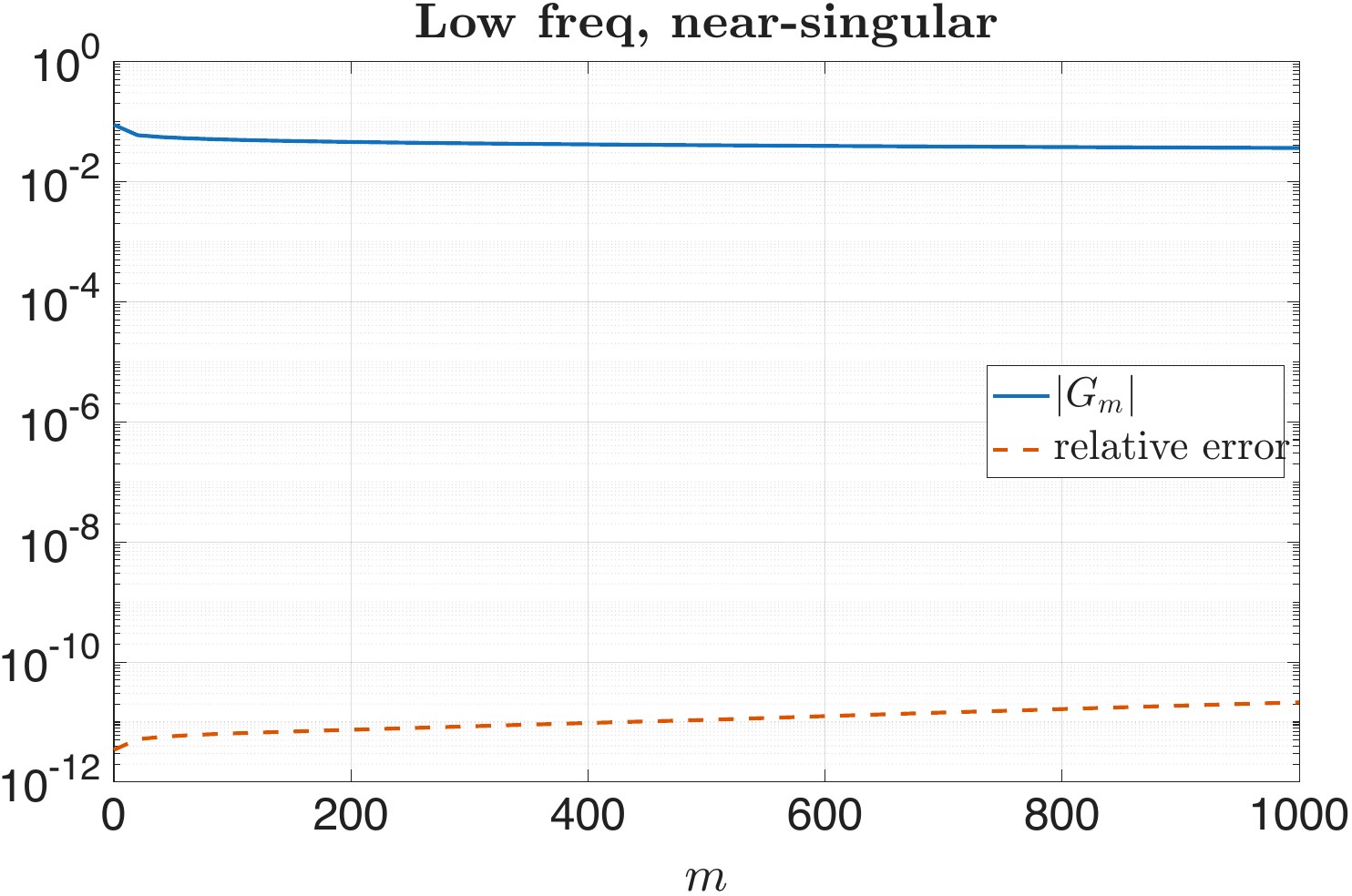}
\caption{Low frequency, near-singular. $\kappa = 6.2 \times 10^{-12}$, $\alpha = 1 - 2.7 \times 10^{-12}$, condition number $= 6.2 \times 10^5$.}
\label{fig:lowfreq_singular}
\end{subfigure}\hfill
\begin{subfigure}[b]{0.48\textwidth}
\includegraphics[width=\textwidth]{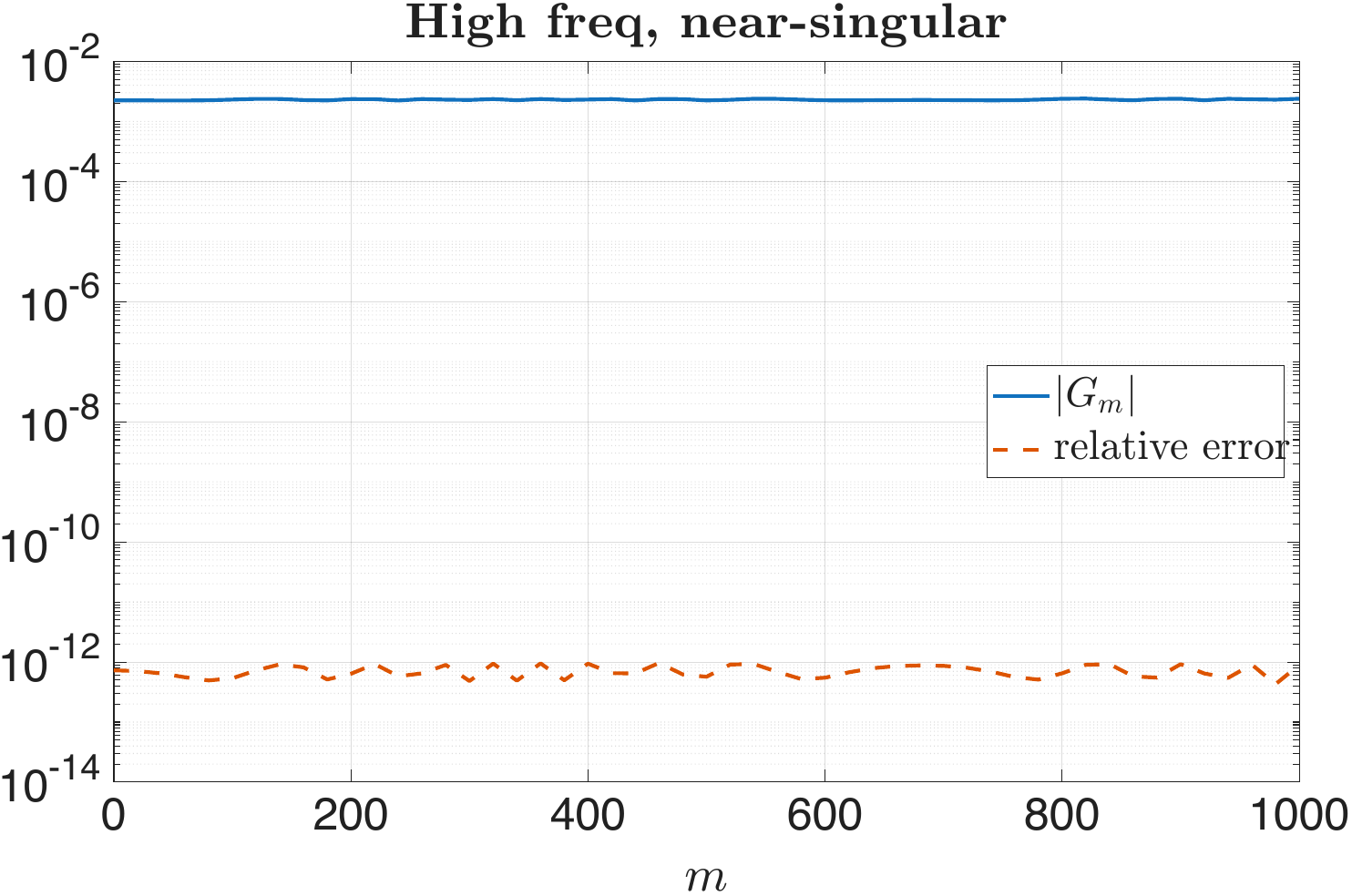}
\caption{High frequency, near-singular. $\kappa = 6152$, $\alpha = 1 - 2.7 \times 10^{-12}$, condition number $= 4.6 \times 10^9$.}
\label{fig:highfreq_singular}
\end{subfigure}

\vspace{0.5em}
\begin{center}
\begin{subfigure}[b]{0.48\textwidth}
\includegraphics[width=\textwidth]{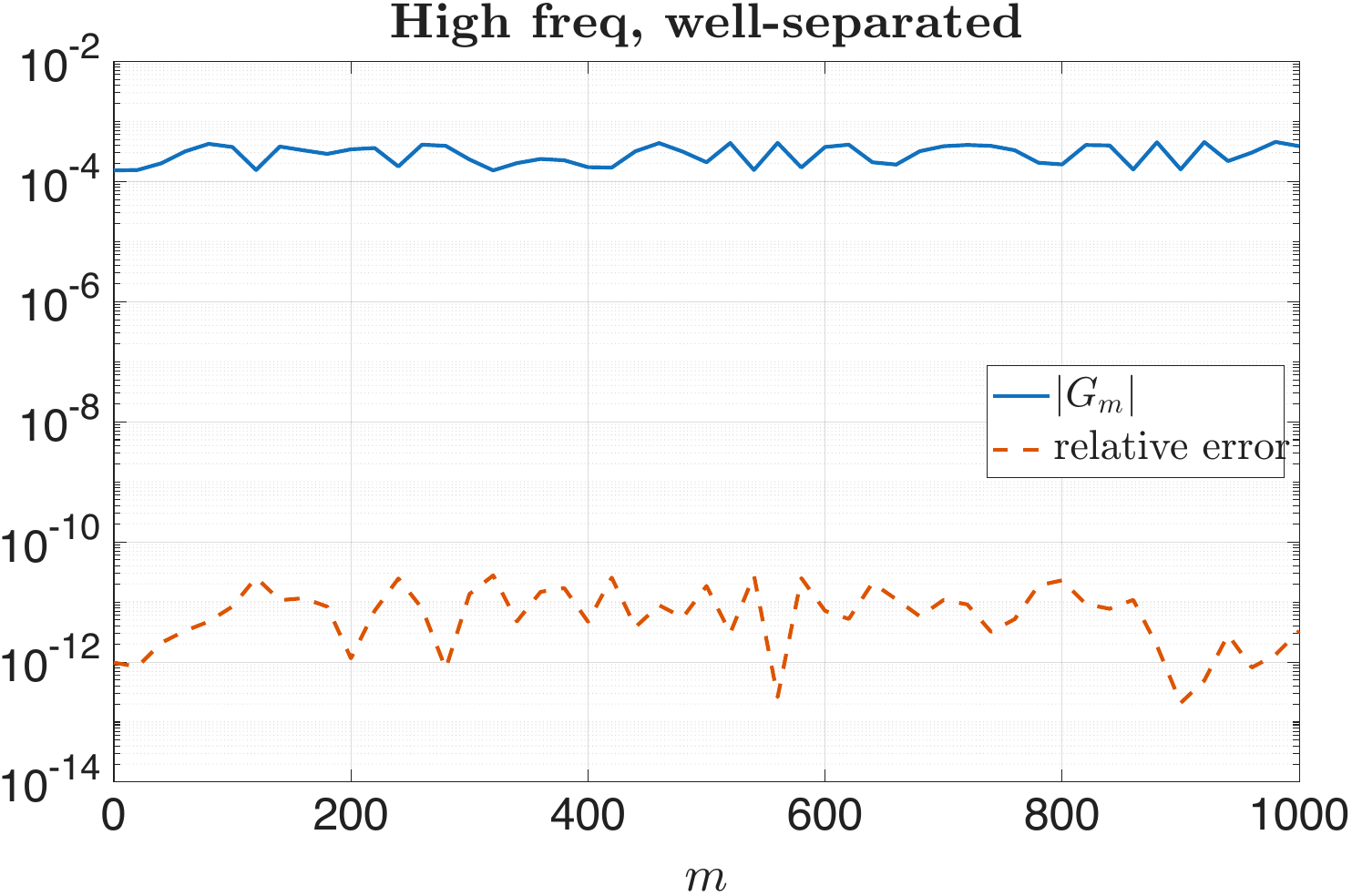}
\caption{High frequency, well-separated. $\kappa = 4380$, $\alpha = 0.902$, condition number $= 1.4 \times 10^9$.}
\label{fig:highfreq_wellsep}
\end{subfigure}
\end{center}
\caption{Accuracy of the pentadiagonal solve. Solid: $|G_m|$; dashed: relative error against adaptive integration in quadruple precision.}
\label{fig:accuracy}
\end{figure}

\begin{figure}[!ht]
\centering
\begin{subfigure}[b]{0.48\textwidth}
\includegraphics[width=\textwidth]{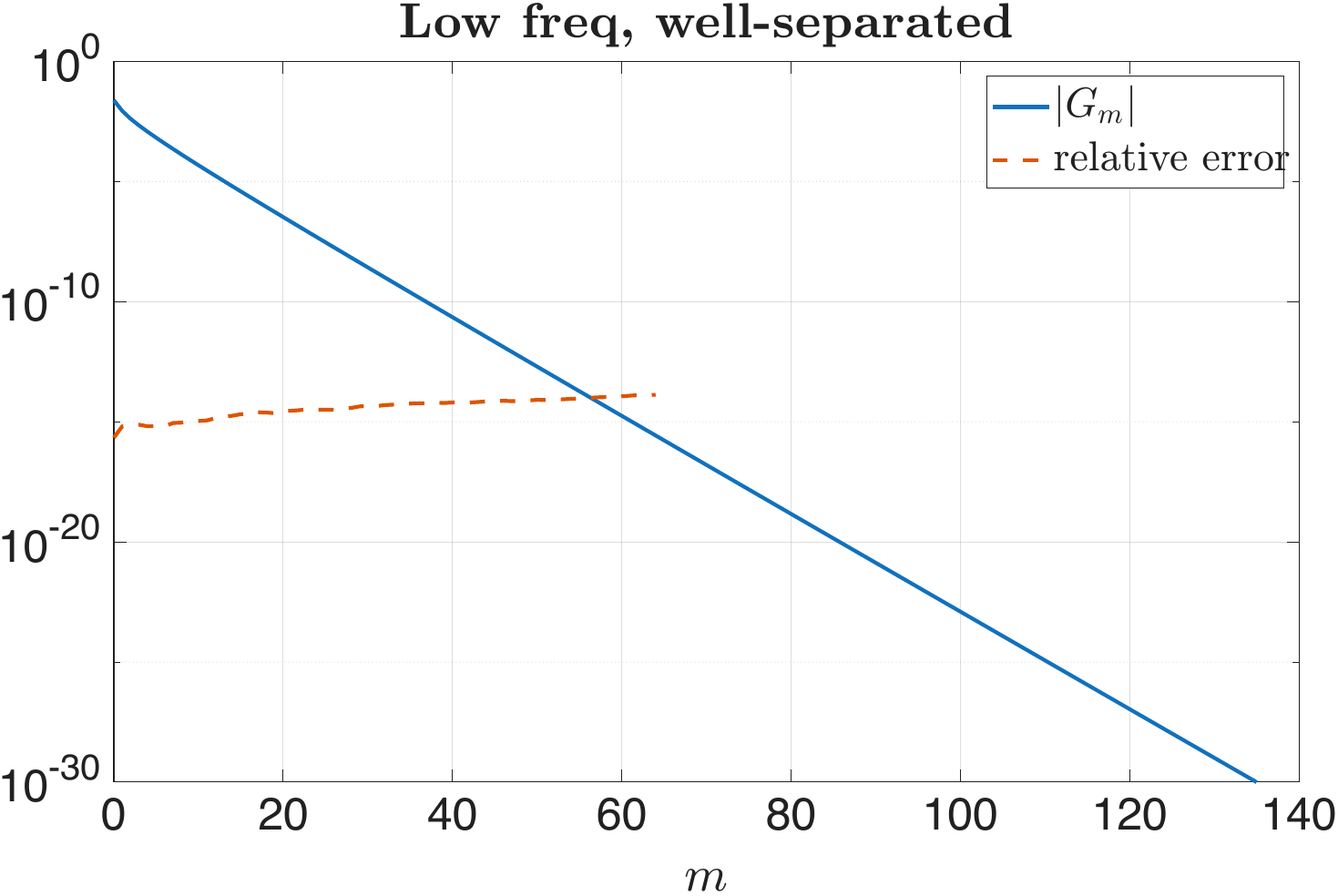}
\caption{Solution obtained with $G_{M'-1} = G_{M'} = 0$.}
\label{fig:lowfreq_wellsep}
\end{subfigure}\hfill
\begin{subfigure}[b]{0.48\textwidth}
\includegraphics[width=\textwidth]{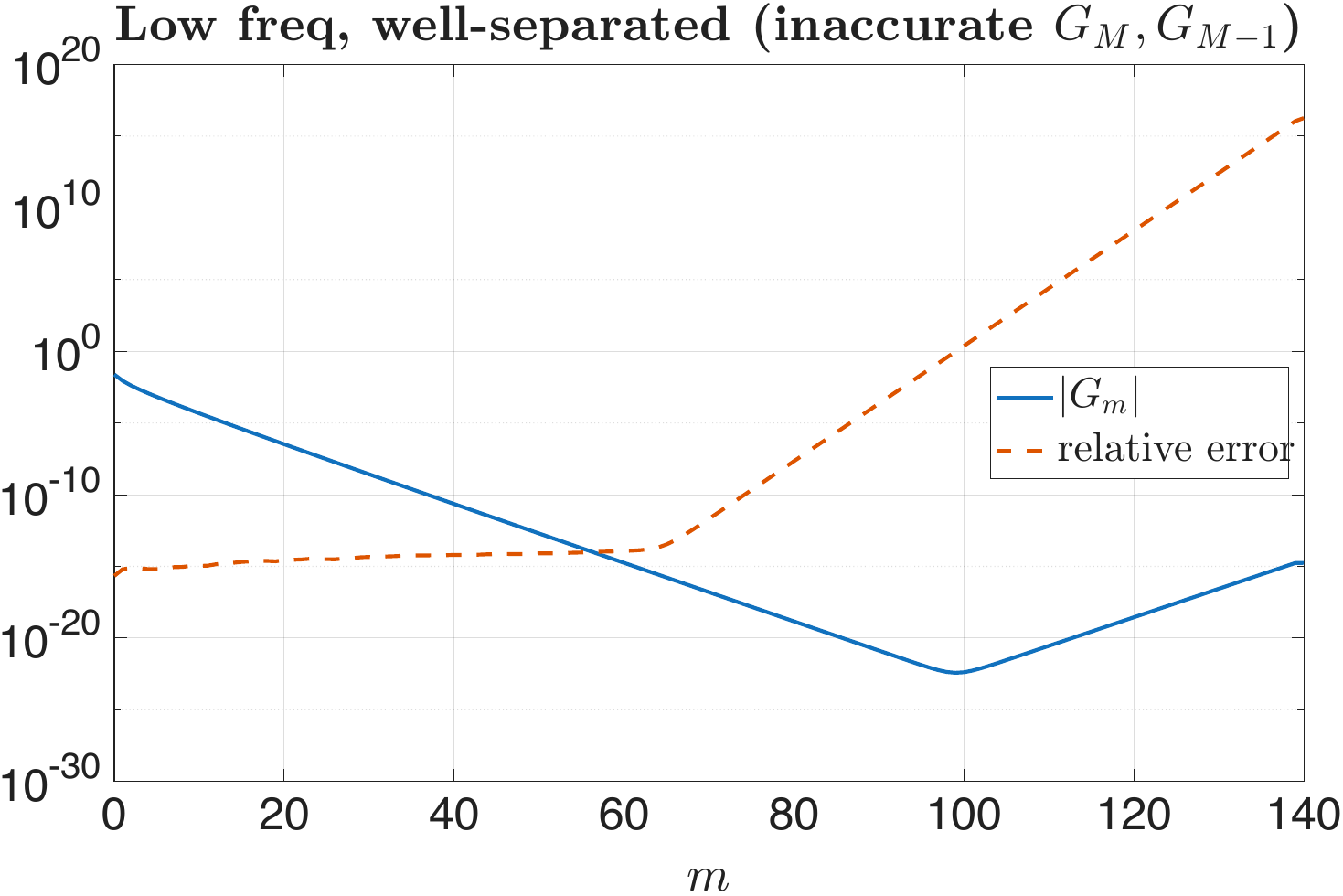}
\caption{Solution obtained with $G_{M'-1}\,, G_{M'}$ from contour deformation.}
\label{fig:lowfreq_wellsep_forced}
\end{subfigure}
\caption{Low frequency, well-separated ($\kappa = 0.438$, $\alpha = 0.902$, cond.\ num.\ $= 19.4$). Solid: $|G_m|$; dashed: relative error against adaptive integration in quadruple precision. In (b), the solid line is the computed $|G_m|$: beyond $m \approx 100$ it no longer follows the exact $|G_m|$, shown in (a), but the error carried in from the inaccurate boundary values.}
\label{fig:wellsep}
\end{figure}

We next test the componentwise accuracy of the pentadiagonal solve in four representative regimes: low and high frequency, each with near-singular and well-separated geometries. In each case we compute $G_m$\,, $m = 0, \ldots, M$\,, with $M=1000$\,, and compare against adaptive integration in quadruple precision. Figures~\ref{fig:accuracy}--\ref{fig:wellsep} show the mode magnitudes and relative errors; Table~\ref{tab:accuracy} reports the right-hand-side error, the condition number of the pentadiagonal matrix, and the maximum relative error over all computed modes.

In all three cases of Figure~\ref{fig:accuracy}, the relative error remains roughly uniform across modes $m = 0, \ldots, M$ near the level of the right-hand-side error, never exceeding it by more than a factor of $\sim 10$, despite condition numbers as large as $10^9$. This reflects the structure of the perturbations: the condition number measures worst-case amplification over arbitrary perturbations, whereas the perturbations arising here are highly structured, since the right-hand side is generated by boundary values of the same recurrence.

\begin{table}[!ht]
\centering
\begin{tabular}{lcccc}
\toprule
Case & RHS error & cond.\ num. & max rel.\ error \\
\midrule
Low freq, singular   & $2.1 \times 10^{-11}$ & $6.2 \times 10^5$ & $2.1 \times 10^{-11}$ \\
Low freq, well-sep   & $6.5 \times 10^{-16}$ & $19.4$ & $1.2 \times 10^{-14}$ \\
High freq, singular  & $8.9 \times 10^{-13}$ & $4.6 \times 10^9$ & $9.6 \times 10^{-13}$ \\
High freq, well-sep  & $3.3 \times 10^{-12}$ & $1.4 \times 10^9$ & $3.0 \times 10^{-11}$ \\
\bottomrule
\end{tabular}
\caption{Accuracy of the pentadiagonal solve for $M = 1000$. The right-hand side (RHS) error is the maximum relative error of the modes $G_0$, $G_1$, $G_{M-1}$, $G_M$ that enter the right-hand side of (\ref{eq:pentadiag}). In the low-frequency well-separated case, the decay-regime algorithm is used and $G_{M-1}$, $G_M$ are not needed; the RHS error reflects only $G_0$, $G_1$.}
\label{tab:accuracy}
\end{table}

Figure~\ref{fig:wellsep} isolates the decay-regime behavior. The geometric decay of $|G_m|$ is visible as $m$ increases, and we choose $M' = 145$\,, for which $|G_{M'}| < 10^{-32}$\,.
In Figure~\ref{fig:wellsep}(a), we set $G_{M'-1}$ and $G_{M'}$ to zero. The resulting solution retains full double-precision relative accuracy as long as the computed value is above machine epsilon and the quadruple-precision reference is itself accurate to double precision. Below this range, the relative error is not plotted.
Figure~\ref{fig:wellsep}(b) shows the failure mode that this Miller-type strategy avoids. If $G_{M'-1}$ and $G_{M'}$ are instead obtained by contour deformation after the modes have decayed far below machine epsilon, they carry no relative accuracy. The pentadiagonal solve then propagates this inaccurate boundary data into the tail, and high relative accuracy is lost for $m \gtrsim 65$\,. The boundary values returned by contour deformation carry an absolute error far larger than $|G_{M'}| < 10^{-32}$\,, and the pentadiagonal solve carries this error into the interior, where it decreases toward smaller $m$\,. For $m \lesssim 65$ it is negligible compared with $|G_m|$\,; for $m \gtrsim 100$ it exceeds $|G_m|$\,, so the computed values -- still all below $10^{-14}$ -- are dominated by the propagated boundary error rather than by the true $G_m$\,. The kink in the solid curve of Figure~\ref{fig:wellsep}(b) marks this crossover; the exact $|G_m|$ continues to decay as in Figure~\ref{fig:wellsep}(a).

Sections~\ref{sec:contour}--\ref{sec:allmodes} together provide an $O(M)$ framework for evaluating $G_m$ and its derivatives in the non-decay and decay regimes. Contour deformation provides a small number of boundary modes, and the pentadiagonal solve recovers the rest with componentwise relative accuracy.

\begin{remark}
	The pentadiagonal system for $\frac{\partial G_m}{\partial a}$ arising from (\ref{eq:fiveterm_da}) has the same matrix as (\ref{eq:pentadiag}), so once $G_m$ is known the same factorization can be reused with the new right-hand side, and the accuracy behavior is essentially the same as for $G_m$.
\end{remark}

%In summary, the pentadiagonal solve preserves the relative accuracy of the boundary values across all modes, even when the matrix condition number exceeds $10^9$ — a componentwise behavior far better than the worst-case bound suggests. 

\section{Algorithm}
\label{sec:algorithm}

We now combine the components of Sections~\ref{sec:contour}--\ref{sec:allmodes} into a complete algorithm. Given a real wavenumber $k \ge 0$\,, a source $(r', z')$ and target $(r, z)$ with $r, r' > 0$\,, and a maximum Fourier mode $M$\,, the algorithm computes $G_m(r,z;r',z')$ for $m = 0, 1, \ldots, M$ (negative modes follow from the symmetry $G_{-m} = G_m$). The algorithm also computes all first- and second-order derivatives with respect to $(r, z, r', z')$\,.

As described in Section~\ref{sec:regimes}, the evaluation divides into three regimes. The near-axis regime is handled by precomputed Taylor-expansion tables. We therefore focus on the non-decay regime (Sections~\ref{sec:allgm}--\ref{sec:allderivatives}) and the decay regime (Section~\ref{sec:decay}). In all cases, the total cost is $O(M)$\,, independent of the wavenumber $k$ and of the source-target separation.

\subsection{Computing all $G_m$}
\label{sec:allgm}

In the non-decay regime, we evaluate $G_m$ at the four boundary modes $m \in \{0, 1, M-1, M\}$ by the contour deformation of Section~\ref{sec:contour}. The intermediate modes $G_2, \ldots, G_{M-2}$ are then recovered by solving the pentadiagonal system (\ref{eq:pentadiag}) of Section~\ref{sec:allmodes}. The low boundary modes $G_0, G_1$ require only a fixed number of quadrature nodes, the high boundary modes $G_{M-1}, G_M$ require $O(M)$ quadrature nodes, and the banded QR solve is also $O(M)$, so the total cost is $O(M)$.

\begin{algorithm}[!ht]
\caption{Computing $G_m$ for all $m = 0, 1, \ldots, M$.}
\label{alg:allgm}
\textbf{Input:} wavenumber $k \ge 0$; source $(r', z')$ and target $(r, z)$ with $r, r' > 0$; maximum mode $M \ge 0$.\\
\textbf{Output:} $\{G_m\}_{m=0}^{M}$.
\begin{algorithmic}[1]
\If{$M \le 5$}
\State Evaluate $G_0, G_1, \ldots, G_M$ directly by contour deformation.
\ElsIf{$M > 5$}
\State Evaluate $G_0$ and $G_1$ via contour deformation.
\State Evaluate $G_{M-1}$ and $G_M$ via contour deformation.
\State Solve the pentadiagonal system (\ref{eq:pentadiag}) for $G_2, \ldots, G_{M-2}$ by banded QR factorization.
\EndIf
\end{algorithmic}
\end{algorithm}

\begin{remark}
\label{rem:seed_pairs}
Several modes can be evaluated on the same contour at essentially the cost of one mode: the contour parametrization, square roots, and exponential are shared across modes, while the Chebyshev polynomials $T_m$ are generated by a single three-term recurrence. We therefore compute $G_0, G_1$ together on the $m=5$ contour (Remark~\ref{rem:lowmode}), and $G_{M-1}, G_M$ together on the $m=M$ contour. When $M \le 5$, all modes are computed on the $m=5$ contour.
\end{remark}

\subsection{Computing all $G_m$ with derivatives}
\label{sec:allderivatives}

We extend Algorithm~\ref{alg:allgm} to compute the first- and second-order derivatives of $G_m$ with respect to $(r, z, r', z')$. The key point is that only $G_m$ itself is recovered from the pentadiagonal boundary-value solve; once $G_m$ is known for $m = 0, 1, \ldots, M$, the derivative quantities are obtained from the recurrences of Section~\ref{sec:recurrence}.

Specifically, the four direct quantities (\ref{eq:directquants}) are evaluated by upward recurrence in $m$ via (\ref{eq:dgda_recurrence}), (\ref{eq:dgdapb_recurrence}), (\ref{eq:d2gda2_recurrence}), and (\ref{eq:d2gda2pab_recurrence}), starting from their values at $m = 0, 1$. These starting values are evaluated by contour deformation (Appendix~\ref{sec:appendix_contour}) alongside $G_0, G_1$ themselves, at little additional quadrature cost (Remark~\ref{rem:shared_cost}). The three remaining quantities in (\ref{eq:recquants}) are then recovered from (\ref{eq:dgdb_recurrence}), (\ref{eq:d2gdadb_recurrence}), and (\ref{eq:d2gdbpbb_recurrence}), and the cylindrical-coordinate derivatives follow from the chain rule of Appendix~\ref{sec:appendix_derivatives}. Since all recurrence steps require only $O(M)$ arithmetic, the total cost matches Algorithm~\ref{alg:allgm}.

\begin{algorithm}[!ht]
\caption{Computing $G_m$ and its first- and second-order derivatives for $m = 0, 1, \ldots, M$.}
\label{alg:allgm_deriv}
\textbf{Input:} wavenumber $k \ge 0$; source $(r', z')$ and target $(r, z)$ with $r, r' > 0$; maximum mode $M \ge 0$.\\
\textbf{Output:} $G_m$ and its first- and second-order derivatives with respect to $(r, z, r', z')$, for $m = 0, 1, \ldots, M$.
\begin{algorithmic}[1]
\If{$M \le 5$}
\State Evaluate $G_m$ for $m = 0, 1, \ldots, M$, and $\frac{\partial G_m}{\partial a}$, $\frac{\partial G_m}{\partial a}+\frac{\partial G_m}{\partial b}$, $\frac{\partial^2 G_m}{\partial a^2}$, $\frac{\partial^2 G_m}{\partial a^2}+\frac{\partial^2 G_m}{\partial a\partial b}$ for $m = 0, 1, \ldots, M+1$ by contour deformation.
\ElsIf{$M > 5$}
\State Evaluate $G_m$, $\frac{\partial G_m}{\partial a}$, $\frac{\partial G_m}{\partial a}+\frac{\partial G_m}{\partial b}$, $\frac{\partial^2 G_m}{\partial a^2}$, and $\frac{\partial^2 G_m}{\partial a^2}+\frac{\partial^2 G_m}{\partial a\partial b}$ at $m = 0, 1$ by contour deformation.
\State Evaluate $G_{M-1}$ and $G_M$ by contour deformation.
\State Solve the pentadiagonal system (\ref{eq:pentadiag}) for $G_2, \ldots, G_{M-2}$ by banded QR factorization.
\State Apply (\ref{eq:dgda_recurrence}) upward to obtain $\frac{\partial G_m}{\partial a}$ for $m = 2, \ldots, M+1$, and (\ref{eq:dgdapb_recurrence}) to obtain $\frac{\partial G_m}{\partial a}+\frac{\partial G_m}{\partial b}$ for $m = 2, \ldots, M$.
\State Apply (\ref{eq:d2gda2_recurrence}) and (\ref{eq:d2gda2pab_recurrence}) upward to obtain $\frac{\partial^2 G_m}{\partial a^2}$ and $\frac{\partial^2 G_m}{\partial a^2}+\frac{\partial^2 G_m}{\partial a\partial b}$ for $m = 2, \ldots, M+1$.
\EndIf
\State Recover $\frac{\partial G_m}{\partial b}$, $\frac{\partial^2 G_m}{\partial a\partial b}$, and $\frac{\partial^2 G_m}{\partial a\partial b}+\frac{\partial^2 G_m}{\partial b^2}$ from (\ref{eq:dgdb_recurrence}), (\ref{eq:d2gdadb_recurrence}), (\ref{eq:d2gdbpbb_recurrence}).
\State Convert all $(a,b)$-derivatives to $(r,z,r',z')$-derivatives via the chain rule formulas of Appendix~\ref{sec:appendix_derivatives}.
\end{algorithmic}
\end{algorithm}

\begin{remark}
If only first-order derivatives are needed, the second-order recurrences (\ref{eq:d2gda2_recurrence})--(\ref{eq:d2gda2pab_recurrence}) and the corresponding pointwise identity (\ref{eq:d2gdbpbb_recurrence}) are skipped, and the contour deformation at $m = 0, 1$ only needs to produce $\partial G_m/\partial a$ and $\partial G_m/\partial a + \partial G_m/\partial b$.
\end{remark}

\begin{remark}
\label{rem:pentadiag_da}
When $M$ exceeds a few hundred, the second-order derivatives computed by Algorithm~\ref{alg:allgm_deriv} gradually lose accuracy. This is because the upward recurrence (\ref{eq:d2gda2_recurrence}) has $\partial G_m/\partial a$ on its right-hand side, and the error in $\partial G_m/\partial a$ grows linearly with $m$ due to the factor $\frac{2m}{b}$ in (\ref{eq:dgda_recurrence}), as is evident from (\ref{eq:twoterm_solution}); the second pass through (\ref{eq:d2gda2_recurrence}) then amplifies this to quadratic growth. For $M > 300$, we mitigate this by also evaluating $\partial G_{M-1}/\partial a$ and $\partial G_M/\partial a$ by contour deformation, and replacing the upward recurrence for $\partial G_m/\partial a$ with a pentadiagonal solve of (\ref{eq:fiveterm_da}), using the same matrix as (\ref{eq:pentadiag}). This largely suppresses the $m$-dependent error growth. The additional cost comes from the contour evaluation of $\partial G_{M-1}/\partial a$ and $\partial G_M/\partial a$, which requires $O(M)$ quadrature nodes, and the additional pentadiagonal solve; together they increase the cost by approximately $30\%$\,.
\end{remark}

\subsection{Decay regime}
\label{sec:decay}

As described in Section~\ref{sec:regimes}, when $M$ exceeds the transition mode
\begin{align*}
m^* = \frac{\kappa}{\sqrt{2}}\sqrt{1 - \sqrt{1-\alpha^2}}\,,
\end{align*}
the modal Green's function $|G_m|$ decays exponentially for $m \ge m^*$. In this regime, evaluating the boundary values $G_{M-1}$ and $G_M$ by contour deformation produces values below machine epsilon, so the pentadiagonal solve loses relative accuracy (see Figure~\ref{fig:wellsep}(b)).

As introduced in Section~\ref{sec:allmodes}, we use a Miller-type strategy \cite{olver1967numerical}: choose $M'$ such that $|G_{M'-1}|$ and $|G_{M'}|$ are well below machine epsilon, set them to zero, and solve the pentadiagonal system of size $(M'-3) \times (M'-3)$ for $G_2, \ldots, G_{M'-2}$. The choice of $M'$ depends on the decay rate, not on $M$, so $M'$ may be smaller or larger than $M$. $G_2, \ldots, G_{M'-2}$ are obtained from the linear solve, $G_0, G_1$ from contour deformation, and $G_{M'-1}, G_{M'}$ are zero by construction. Any remaining modes $G_m$ with $m > M'$ are taken to be zero.

In the decay regime, the upward recurrence loses relative accuracy beyond the transition point $m^*$, so a downward recurrence is used instead. Specifically, for each recurrence of Section~\ref{sec:allderivatives}, we set the values at $m = M', M'+1$ to zero and evaluate downward to $m = 0$. Algorithms~\ref{alg:decay_gm} and~\ref{alg:decay_deriv} are the decay-regime counterparts to Algorithms~\ref{alg:allgm} and~\ref{alg:allgm_deriv}, respectively.

\begin{algorithm}[!ht]
\caption{Decay-regime counterpart to Algorithm~\ref{alg:allgm}: computing $G_m$ for $m = 0, 1, \ldots, M$ when $M > m^*$.}
\label{alg:decay_gm}
\textbf{Input:} same as Algorithm~\ref{alg:allgm}, with $M > m^*$.\\
\textbf{Output:} $\{G_m\}_{m=0}^{M}$.
\begin{algorithmic}[1]
\State Choose $M'$ such that $|G_{M'-1}|$ and $|G_{M'}|$ are well below machine epsilon (estimated from $m^*$).
\State Evaluate $G_0$ and $G_1$ by contour deformation.
\State Set $G_{M'-1} = G_{M'} = 0$ and solve the pentadiagonal system (\ref{eq:pentadiag}) of size $(M'-3) \times (M'-3)$ for $G_2, \ldots, G_{M'-2}$ by banded QR factorization.
\State If $M > M'$, set $G_m = 0$ for $m = M'+1, \ldots, M$.
\end{algorithmic}
\end{algorithm}

\begin{algorithm}[!ht]
\caption{Decay-regime counterpart to Algorithm~\ref{alg:allgm_deriv}: computing $G_m$ and its first- and second-order derivatives for $m = 0, 1, \ldots, M$ when $M > m^*$.}
\label{alg:decay_deriv}
\textbf{Input:} same as Algorithm~\ref{alg:allgm_deriv}, with $M > m^*$.\\
\textbf{Output:} $G_m$ and its first- and second-order derivatives with respect to $(r, z, r', z')$, for $m = 0, 1, \ldots, M$.
\begin{algorithmic}[1]
\State Choose $M'$ such that $|G_{M'-1}|$ and $|G_{M'}|$ are well below machine epsilon (estimated from $m^*$).
\State Evaluate $G_0$ and $G_1$ by contour deformation.
\State Set $G_{M'-1} = G_{M'} = 0$ and solve the pentadiagonal system (\ref{eq:pentadiag}) of size $(M'-3) \times (M'-3)$ for $G_2, \ldots, G_{M'-2}$ by banded QR factorization.
\State Set each derivative quantity at $m = M', M'+1$ to zero. Apply (\ref{eq:dgda_recurrence}) and (\ref{eq:dgdapb_recurrence}) downward to obtain $\frac{\partial G_m}{\partial a}$ and $\frac{\partial G_m}{\partial a}+\frac{\partial G_m}{\partial b}$ for $m = 0, \ldots, M'-1$.
\State Apply (\ref{eq:d2gda2_recurrence}) and (\ref{eq:d2gda2pab_recurrence}) downward to obtain $\frac{\partial^2 G_m}{\partial a^2}$ and $\frac{\partial^2 G_m}{\partial a^2}+\frac{\partial^2 G_m}{\partial a\partial b}$ for $m = 0, \ldots, M'-1$.
\State Recover $\frac{\partial G_m}{\partial b}$, $\frac{\partial^2 G_m}{\partial a\partial b}$, and $\frac{\partial^2 G_m}{\partial a\partial b}+\frac{\partial^2 G_m}{\partial b^2}$ from (\ref{eq:dgdb_recurrence}), (\ref{eq:d2gdadb_recurrence}), (\ref{eq:d2gdbpbb_recurrence}).
\State If $M > M'$, set $G_m$ and all its derivatives to zero for $m = M'+1, \ldots, M$.
\State Convert all $(a,b)$-derivatives to $(r,z,r',z')$-derivatives via the chain rule formulas of Appendix~\ref{sec:appendix_derivatives}.
\end{algorithmic}
\end{algorithm}

\section{Numerical results}
\label{sec:numerical}

We present numerical experiments for the modal Green's function evaluator and for its use in boundary integral equation solvers for acoustic scattering from bodies of revolution. Section~\ref{sec:accuracy_timing} tests the evaluator on its own, including accuracy in the non-decay and decay regimes, independence of the wavenumber and source-target separation, and the observed $O(M)$ scaling. Section~\ref{sec:bie} then incorporates the evaluator into modal BIE solvers and measures the resulting assembly and solve costs. All experiments are performed in Fortran~77, compiled with \texttt{gfortran -O3}, on a 2021 MacBook Pro (Apple M1 Max, 64\,GB RAM). Timings in Section~\ref{sec:accuracy_timing} are single-core; Section~\ref{sec:bie} uses up to 8 cores.

\subsection{Accuracy and timing}
\label{sec:accuracy_timing}

\begin{figure}[!ht]
\centering
\begin{subfigure}[b]{0.48\textwidth}
\includegraphics[width=\textwidth]{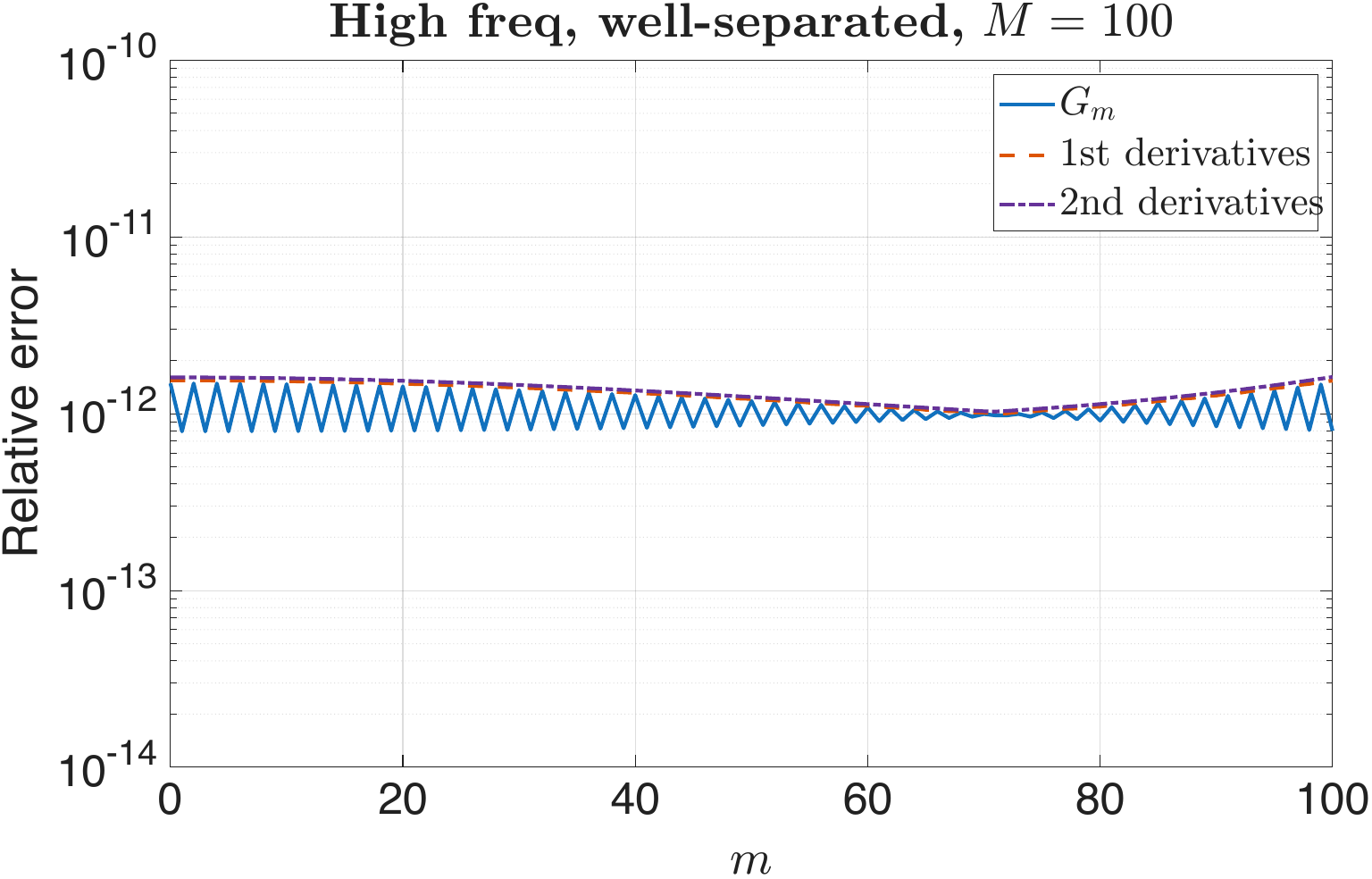}
\caption{Well-separated, $M = 100$.}
\label{fig:deriv_wellsep_M100}
\end{subfigure}\hfill
\begin{subfigure}[b]{0.48\textwidth}
\includegraphics[width=\textwidth]{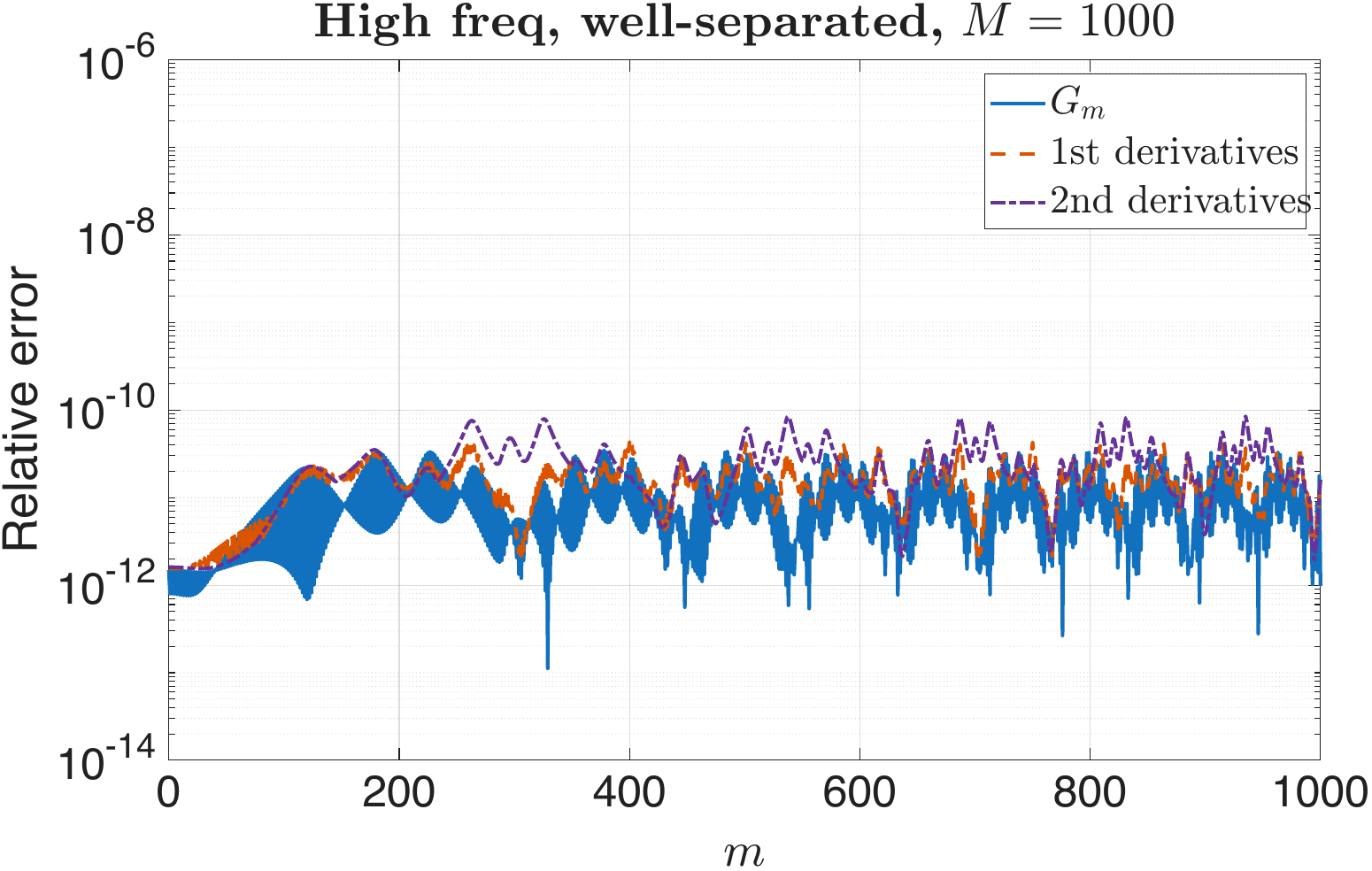}
\caption{Well-separated, $M = 1000$.}
\label{fig:deriv_wellsep_M1000}
\end{subfigure}

\vspace{0.5em}
\begin{subfigure}[b]{0.48\textwidth}
\includegraphics[width=\textwidth]{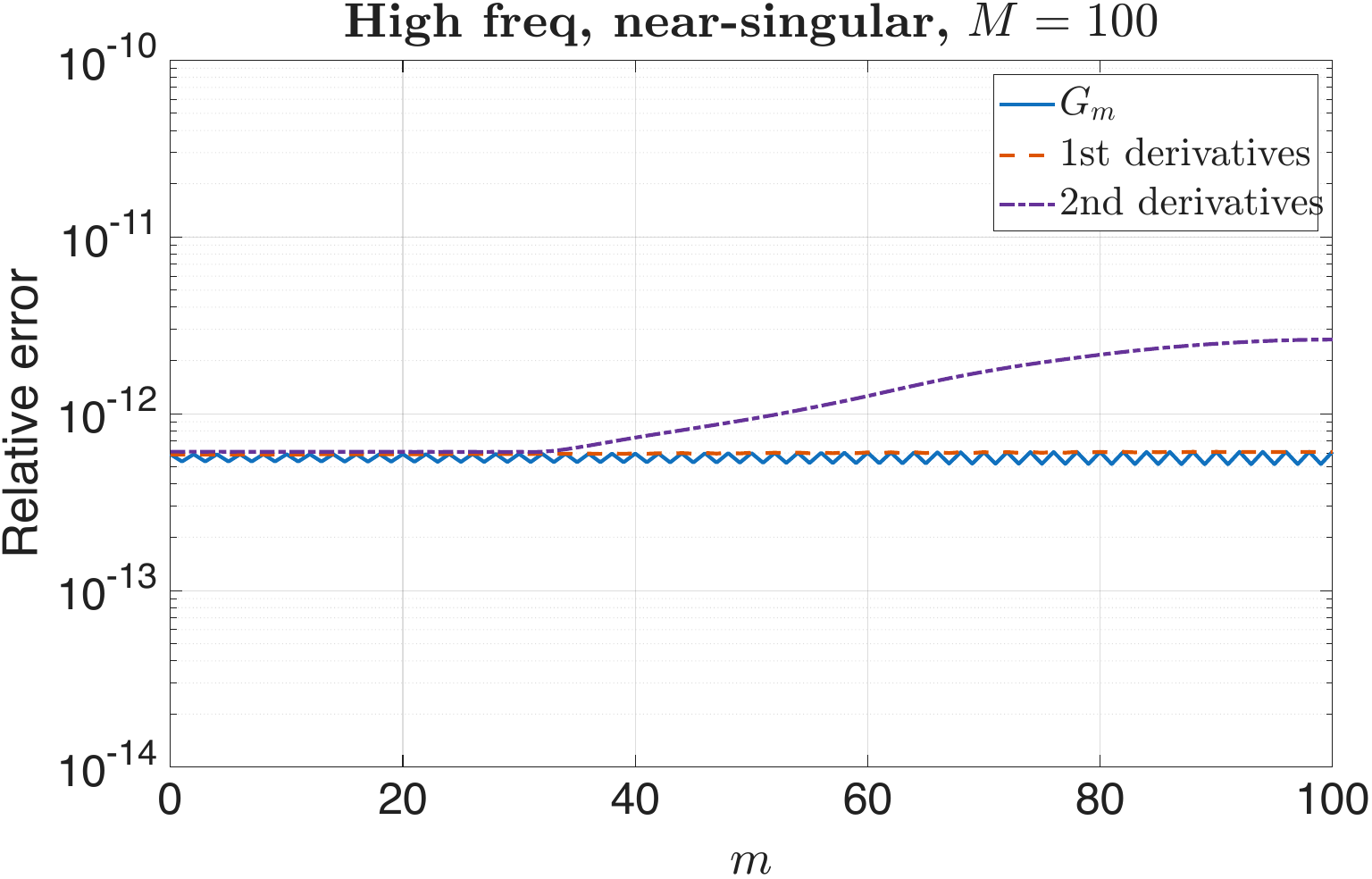}
\caption{Near-singular, $M = 100$.}
\label{fig:deriv_singular_M100}
\end{subfigure}\hfill
\begin{subfigure}[b]{0.48\textwidth}
\includegraphics[width=\textwidth]{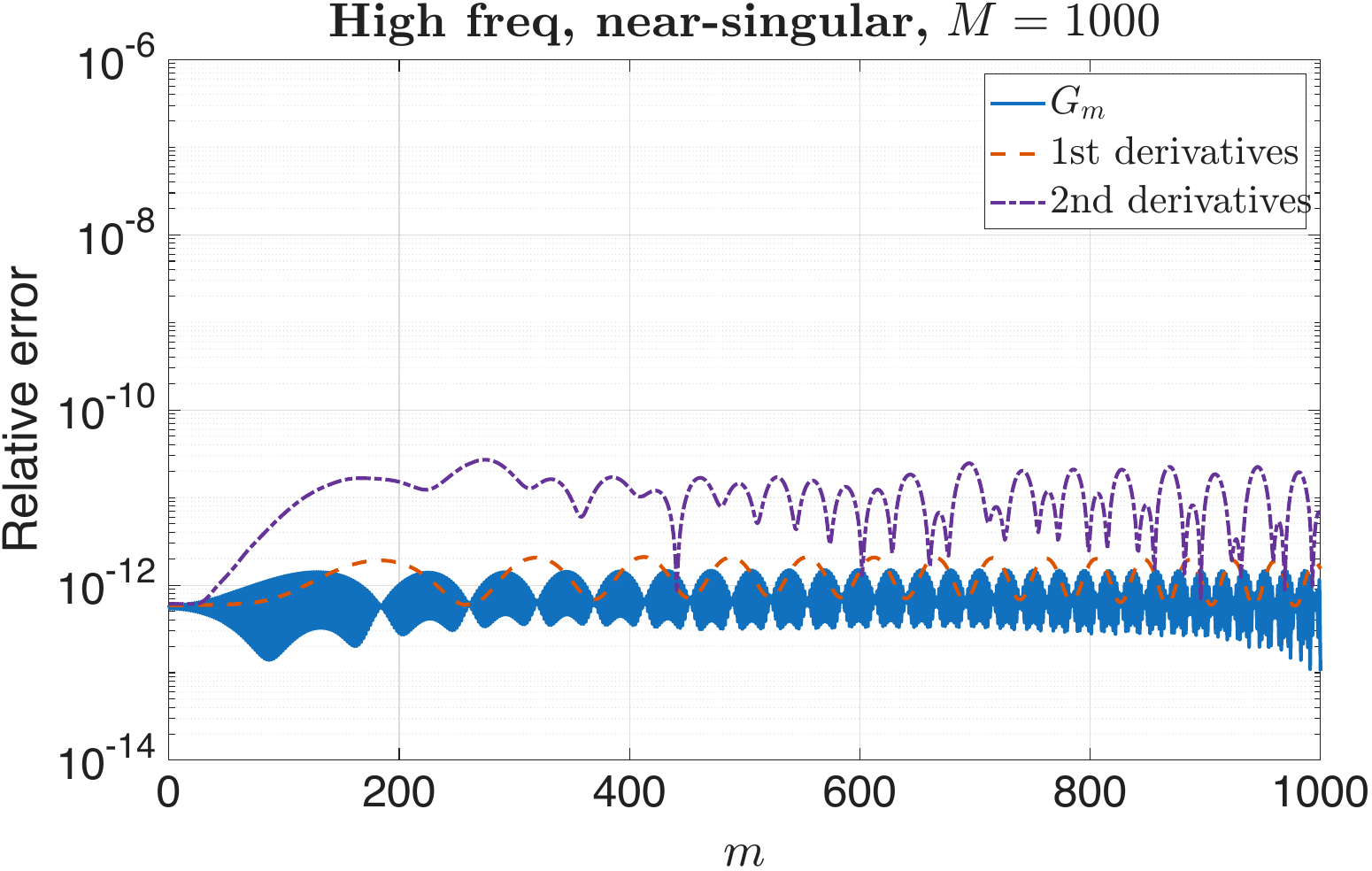}
\caption{Near-singular, $M = 1000$.}
\label{fig:deriv_singular_M1000}
\end{subfigure}
\caption{Non-decay regime accuracy for $k = 2500$\,. Top row: well-separated geometry ($\kappa = 10949$, $\alpha = 0.902$). Bottom row: near-singular geometry ($\kappa = 15397$, $\alpha = 1 - 2.7 \times 10^{-12}$). The curves show relative errors in $G_m$ (solid), first-order derivatives (dashed), and second-order derivatives (dash-dotted).}
\label{fig:deriv_accuracy}
\end{figure}

\begin{figure}[!ht]
\centering
\begin{subfigure}[b]{0.48\textwidth}
\includegraphics[width=\textwidth]{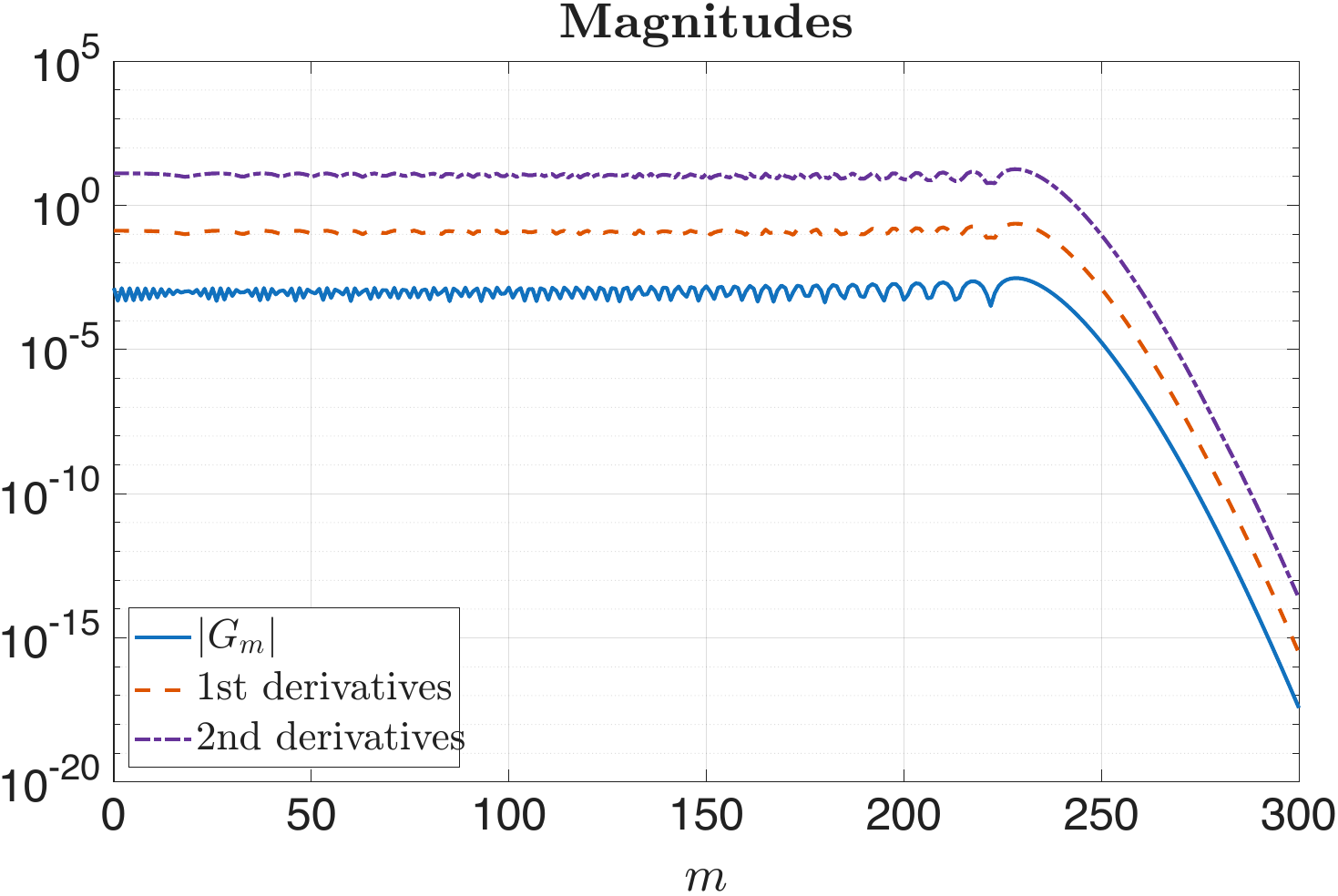}
\caption{Magnitudes.}
\label{fig:decay_mag}
\end{subfigure}\hfill
\begin{subfigure}[b]{0.48\textwidth}
\includegraphics[width=\textwidth]{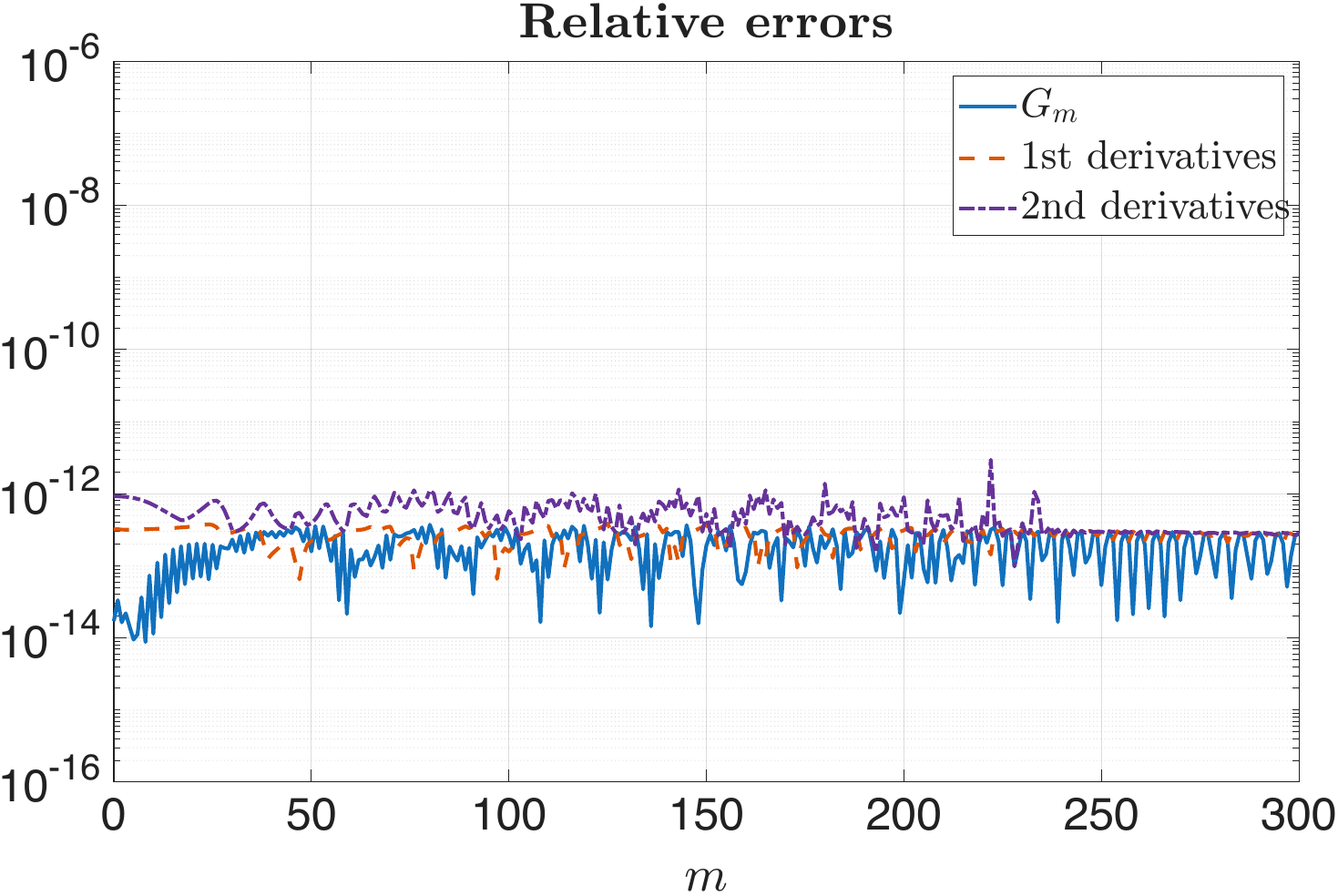}
\caption{Relative errors.}
\label{fig:decay_err}
\end{subfigure}
\caption{Decay-regime accuracy for $k = 100$\,, $\kappa = 438$\,, $\alpha = 0.902$\,, $M = 300$\,, with transition point $m^* = 233$\,. Left: magnitudes of $G_m$\,, first-order derivatives, and second-order derivatives. Right: relative errors. Relative accuracy is retained across more than 15 orders of magnitude of decay.}
\label{fig:decay}
\end{figure}

\begin{figure}[!ht]
\centering
\includegraphics[width=0.55\textwidth]{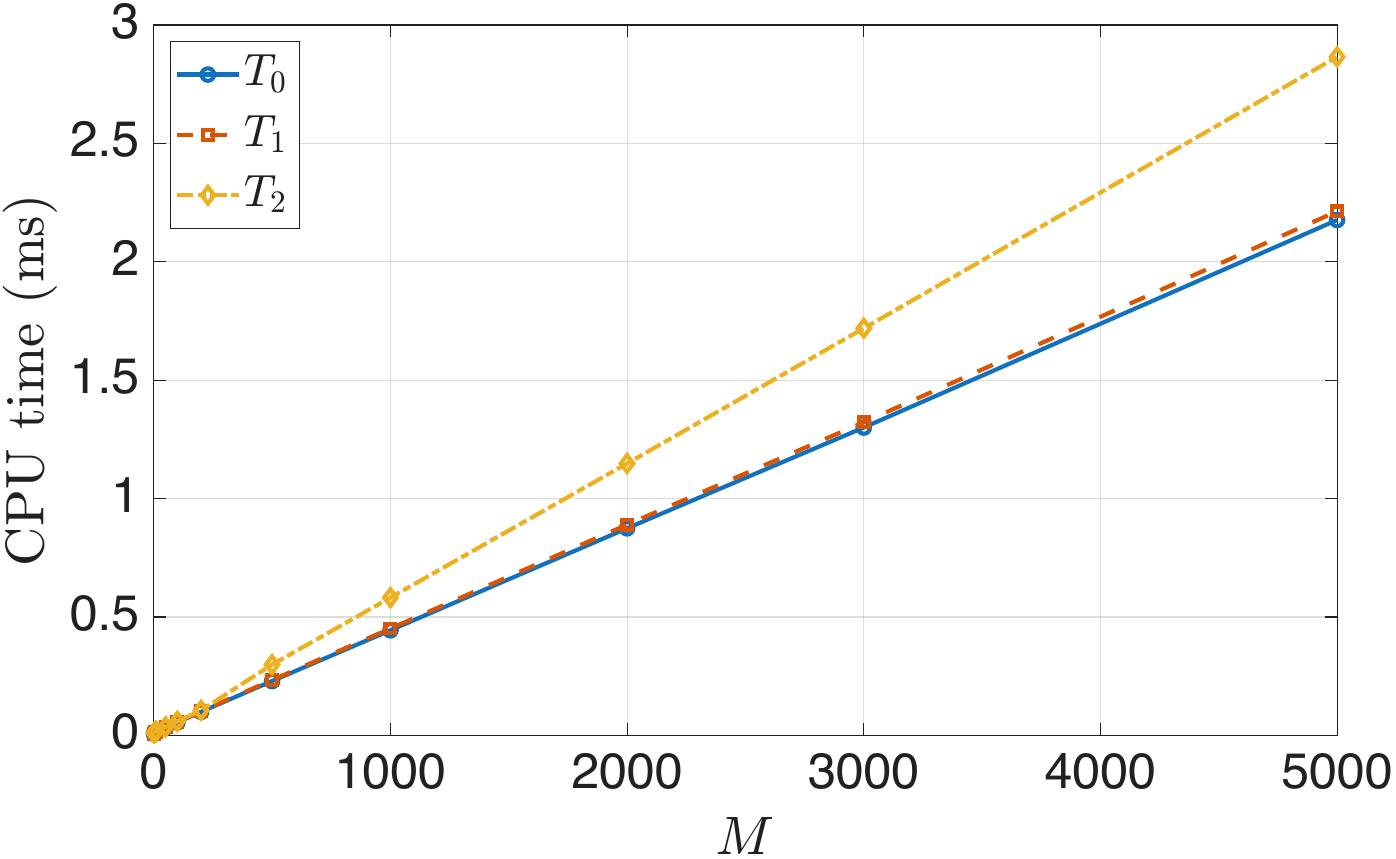}
\caption{Single-core CPU time as a function of $M$ for $k = 2500$\,, $\alpha = 0.902$\,. Here $T_0$ denotes evaluation of all $G_m$\,, $T_1$ includes first-order derivatives, and $T_2$ includes first- and second-order derivatives.}
\label{fig:scaling}
\end{figure}

The experiments are organized as follows. Figure~\ref{fig:deriv_accuracy} and Table~\ref{tab:deriv_accuracy} test accuracy in the non-decay regime for $M$ up to $3000$, in both well-separated and near-singular geometries. Figure~\ref{fig:decay} tests accuracy in the decay regime. Tables~\ref{tab:timing_k}--\ref{tab:timing_geom} measure the dependence of the runtime on $k$ and on the singularity parameter $\alpha$, and Figure~\ref{fig:scaling} verifies the $O(M)$ scaling. All reference values are computed by adaptive integration in quadruple precision.

We begin with the non-decay regime at $k = 2500$\,. Figure~\ref{fig:deriv_accuracy} shows, for each mode $m$\,, the relative error in $G_m$ and its first- and second-order derivatives. For $M = 100$\,, all quantities are accurate to about 12 digits uniformly across all modes, and the derivative errors are essentially the same size as the errors in $G_m$\,. For $M = 1000$\,, the accuracy remains nearly uniform, with $G_m$ and its first- and second-order derivatives all retaining about 11 digits. Even at $M = 3000$\,, all quantities retain about 10--11 digits of relative accuracy (Table~\ref{tab:deriv_accuracy}); the gradual loss as $M$ increases reflects the conditioning of the high-mode boundary-value integrals, whose integrands oscillate $O(M)$ times.

\begin{table}[!ht]
\centering
\begin{tabular}{r|ccc|ccc}
\toprule
& \multicolumn{3}{c|}{Well-separated ($\alpha = 0.902$)} & \multicolumn{3}{c}{Near-singular ($\alpha \approx 1$)} \\
$M$ & $G_m$ & 1st deriv & 2nd deriv & $G_m$ & 1st deriv & 2nd deriv \\
\midrule
100  & $1.5 \times 10^{-12}$ & $1.5 \times 10^{-12}$ & $1.6 \times 10^{-12}$ & $6.1 \times 10^{-13}$ & $6.1 \times 10^{-13}$ & $2.6 \times 10^{-12}$ \\
1000 & $3.5 \times 10^{-11}$ & $4.3 \times 10^{-11}$ & $8.5 \times 10^{-11}$ & $1.5 \times 10^{-12}$ & $2.1 \times 10^{-12}$ & $2.7 \times 10^{-11}$ \\
2000 & $3.7 \times 10^{-11}$ & $6.4 \times 10^{-11}$ & $1.1 \times 10^{-10}$ & $2.5 \times 10^{-12}$ & $2.7 \times 10^{-12}$ & $2.7 \times 10^{-11}$ \\
3000 & $2.3 \times 10^{-11}$ & $4.7 \times 10^{-11}$ & $8.7 \times 10^{-11}$ & $3.0 \times 10^{-12}$ & $4.1 \times 10^{-12}$ & $4.7 \times 10^{-11}$ \\
\bottomrule
\end{tabular}
\caption{Maximum relative error over all modes $m = 0, \ldots, M$ and all derivative components, for $k = 2500$\,. Well-separated: $\kappa = 10949$\,, $\alpha = 0.902$\,. Near-singular: $\kappa = 15397$\,, $\alpha = 1 - 2.7 \times 10^{-12}$\,.}
\label{tab:deriv_accuracy}
\end{table}

Figure~\ref{fig:decay} tests the decay regime with $k = 100$\,, $\kappa = 438$\,, $\alpha = 0.902$\,, and $M = 300$\,. The transition point, given by (\ref{eq:mstar}), is $m^* = 233$, after which $G_m$ and its derivatives decay exponentially. The magnitudes span more than $15$ orders of magnitude, yet the relative errors remain essentially flat, at about $10^{-12}$\,, down to modes with $|G_m| \approx 10^{-15}$. This confirms that the pentadiagonal solve, with carefully chosen boundary values, preserves relative accuracy through the exponentially decaying tail.

For the timing tests, let $T_0$\,, $T_1$\,, and $T_2$ denote the single-core CPU time for evaluating all modes $G_m$\,, $m = 0, \ldots, M$\,, with no derivatives, with first-order derivatives, and with both first- and second-order derivatives, respectively. All timings are averaged over $10{,}000$ evaluations.

Across wavenumbers $k = 10$ to $5000$ at the well-separated geometry ($\alpha = 0.902$) of Figure~\ref{fig:deriv_accuracy}, timings are essentially independent of~$k$ (Table~\ref{tab:timing_k}). At fixed $k = 2500$, sweeping the source-target distance from $d = 1.37$ down to $d = 1.1 \times 10^{-6}$ (more singular than Figure~\ref{fig:deriv_accuracy}, where $\alpha \approx 1 - 2.7 \times 10^{-12}$), timings remain essentially independent of separation, even in the near-singular regime (Table~\ref{tab:timing_geom}).

\begin{table}[!ht]
\centering
\footnotesize
\begin{tabular}{r|ccc|ccc|ccc|ccc}
\toprule
& \multicolumn{3}{c|}{$M = 10$} & \multicolumn{3}{c|}{$M = 100$} & \multicolumn{3}{c|}{$M = 1000$} & \multicolumn{3}{c}{$M = 5000$} \\
$k$ & $T_0$ & $T_1$ & $T_2$ & $T_0$ & $T_1$ & $T_2$ & $T_0$ & $T_1$ & $T_2$ & $T_0$ & $T_1$ & $T_2$ \\
\midrule
10   & 0.015 & 0.016 & 0.017 & 0.056 & 0.055 & 0.057 & 0.43 & 0.44 & 0.57 & 2.10 & 2.13 & 2.79 \\
100  & 0.015 & 0.016 & 0.017 & 0.054 & 0.056 & 0.058 & 0.43 & 0.44 & 0.57 & 2.10 & 2.14 & 2.79 \\
500  & 0.015 & 0.016 & 0.017 & 0.055 & 0.057 & 0.059 & 0.44 & 0.45 & 0.58 & 2.12 & 2.15 & 2.81 \\
1000 & 0.015 & 0.016 & 0.017 & 0.055 & 0.057 & 0.059 & 0.44 & 0.45 & 0.58 & 2.14 & 2.18 & 2.83 \\
2500 & 0.015 & 0.016 & 0.017 & 0.056 & 0.057 & 0.059 & 0.44 & 0.45 & 0.58 & 2.18 & 2.22 & 2.88 \\
5000 & 0.015 & 0.016 & 0.017 & 0.056 & 0.057 & 0.059 & 0.45 & 0.46 & 0.58 & 2.18 & 2.21 & 2.86 \\
\bottomrule
\end{tabular}
\caption{CPU time (ms) as a function of wavenumber~$k$\,. Fixed geometry: $r = 2.35$\,, $r' = 3.68$\,, $z = 3.16$\,, $z' = 2.82$ ($\alpha = 0.902$).}
\label{tab:timing_k}
\end{table}

\begin{table}[!ht]
\centering
\footnotesize
\begin{tabular}{r|ccc|ccc|ccc|ccc}
\toprule
& \multicolumn{3}{c|}{$M = 10$} & \multicolumn{3}{c|}{$M = 100$} & \multicolumn{3}{c|}{$M = 1000$} & \multicolumn{3}{c}{$M = 5000$} \\
$\alpha$ & $T_0$ & $T_1$ & $T_2$ & $T_0$ & $T_1$ & $T_2$ & $T_0$ & $T_1$ & $T_2$ & $T_0$ & $T_1$ & $T_2$ \\
\midrule
0.902 & 0.015 & 0.016 & 0.017 & 0.055 & 0.057 & 0.059 & 0.44 & 0.45 & 0.58 & 2.18 & 2.23 & 2.87 \\
0.975 & 0.015 & 0.016 & 0.017 & 0.055 & 0.057 & 0.059 & 0.44 & 0.45 & 0.58 & 2.18 & 2.23 & 2.86 \\
0.999 & 0.015 & 0.016 & 0.017 & 0.055 & 0.057 & 0.059 & 0.44 & 0.45 & 0.58 & 2.19 & 2.23 & 2.86 \\
$1 - 10^{-5}$ & 0.015 & 0.016 & 0.017 & 0.055 & 0.056 & 0.059 & 0.44 & 0.45 & 0.58 & 2.19 & 2.23 & 2.88 \\
$1 - 10^{-8}$ & 0.015 & 0.016 & 0.017 & 0.060 & 0.062 & 0.064 & 0.45 & 0.46 & 0.59 & 2.19 & 2.24 & 2.88 \\
$1 - 10^{-14}$ & 0.015 & 0.016 & 0.017 & 0.061 & 0.063 & 0.066 & 0.45 & 0.46 & 0.59 & 2.20 & 2.24 & 2.88 \\
\bottomrule
\end{tabular}
\caption{CPU time (ms) as a function of the singularity parameter~$\alpha$\,. Fixed $k = 2500$\,, $r = 2.35$\,, $z = 3.16$\,. The source-target distance ranges from $d = 1.37$ ($\alpha = 0.902$) to $d = 1.1 \times 10^{-6}$ ($\alpha \approx 1$).}
\label{tab:timing_geom}
\end{table}

Finally, Figure~\ref{fig:scaling} confirms the linear $O(M)$ scaling. The $T_0$ and $T_1$ curves nearly coincide, indicating that first-order derivatives add negligible overhead. Including second-order derivatives increases the runtime by about $30\%$ for $M \geq 500$ due to the additional contour evaluations of $\partial G_{M-1}/\partial a$ and $\partial G_M/\partial a$ and the pentadiagonal solve (Remark~\ref{rem:pentadiag_da}).

\subsection{Application to boundary integral equations}
\label{sec:bie}

We apply the modal Green's function evaluator to BIE solvers for acoustic scattering from bodies of revolution at wavenumber~$k$\,. The generating curve is discretized into $N$ points using 16th-order Gauss-Legendre panels with approximately $12$ points per wavelength and log-singular quadrature corrections \cite{kolm2001numerical, hao2014high}. For each Fourier mode $m$, a dense system is assembled using the modal Green's function $G_m$ and its derivatives, and solved by LU factorization (LAPACK). The purpose of these experiments is not to optimize the linear algebra, but to isolate the effect of the kernel evaluator on the cost of dense modal BIE assembly. We therefore deliberately use standard dense per-mode assembly and LAPACK solves. This makes the kernel-evaluation cost visible and provides a conservative baseline for future coupling with accelerated linear algebra.

Section~\ref{sec:bie_full3d} verifies the solver on two BIE formulations: the combined-field integral equation (CFIE) for the exterior Dirichlet problem, which uses first-order derivatives of $G_m$, and M\"uller's formulation for the transmission problem, which uses second-order derivatives. Section~\ref{sec:bie_kscan} then sweeps the wavenumber at fixed $M$ to demonstrate that the dense LU solve asymptotically dominates assembly.

\subsubsection{Modal BIE solves for full 3D scattering}
\label{sec:bie_full3d}

For both formulations, the boundary data are decomposed into Fourier modes by FFT, with the truncation $M$ chosen to resolve the azimuthal direction.

\paragraph{Exterior Dirichlet problem.}
We use the combined-field representation $u = (D_k + \I k\, S_k)\sigma$\,, where
\begin{equation}
\begin{aligned}
	(S_k\sigma)(\bm{x}) &= \int_S G_k(|\bm{x}-\bm{y}|)\,\sigma(\bm{y})\, \D S(\bm{y})\,, \\
	(D_k\sigma)(\bm{x}) &= \int_S \frac{\partial G_k(|\bm{x}-\bm{y}|)}{\partial \bm{n}'}\,\sigma(\bm{y})\, \D S(\bm{y})
\end{aligned}
\label{eq:sk_dk}
\end{equation}
are the single- and double-layer potentials, so that the Dirichlet data $f$ on $S$ determine the density $\sigma$ through
\begin{align}
	\tfrac{1}{2}\sigma + (D_k + \I k\, S_k)\sigma = f \quad \text{on } S\,.
	\label{eq:cfie}
\end{align}
Expanding $\sigma$ and $f$ in Fourier modes as in (\ref{eq:bie_modal}), each mode satisfies
\begin{align}
	\tfrac{1}{2}\sigma_m(r,z) + 2\pi\int_{\gamma} \Big[\frac{\partial G_{k,m}}{\partial \bm{n}'} + \I k\, G_{k,m}\Big](r,z,r',z')\, \sigma_m(r',z')\, r'\,\D\ell' = f_m(r,z)\,.
	\label{eq:cfie_modal}
\end{align}%
For the sound-soft scattering problem, the CFIE requires $G_m$ and its normal derivative $\frac{\partial G_m}{\partial \bm{n}'}$ at each source-target pair, where $\bm{n}'$ denotes the outward normal to the generating curve at the source. The BIE system matrix is $N \times N$ for each mode, and the algorithm of Section~\ref{sec:allderivatives} supplies $G_m$ and its first-order derivatives for all modes simultaneously.

\paragraph{Transmission problem.}
For the acoustic transmission problem (penetrable obstacle), we use M\"uller's integral equation \cite{kressroach1978transmission, vico2014boundary}, which is built from differences of $G_m$\,, $\frac{\partial G_m}{\partial \bm{n}}$\,, $\frac{\partial G_m}{\partial \bm{n}'}$\,, and the hypersingular kernel $\frac{\partial^2 G_m}{\partial \bm{n}\,\partial \bm{n}'}$ evaluated at the exterior and interior wavenumbers $k_{\mathrm{ext}}$ and $k_{\mathrm{int}}$\,, where $\bm{n}$ and $\bm{n}'$ are the outward normals at the target and source. Specifically, let
\begin{equation}
\begin{aligned}
	(D'_k\tau)(\bm{x}) &= \int_S \frac{\partial G_k(|\bm{x}-\bm{y}|)}{\partial \bm{n}}\,\tau(\bm{y})\, \D S(\bm{y})\,, \\
	(T_k\sigma)(\bm{x}) &= \int_S \frac{\partial^2 G_k(|\bm{x}-\bm{y}|)}{\partial \bm{n}\,\partial \bm{n}'}\,\sigma(\bm{y})\, \D S(\bm{y})
\end{aligned}
\label{eq:dpk_tk}
\end{equation}
denote the normal derivative of the single-layer potential and the normal derivative of the double-layer potential, with $S_k$ and $D_k$ as in (\ref{eq:sk_dk}). With the representation $u_{\mathrm{ext}} = D_{k_{\mathrm{ext}}}\sigma + k_{\mathrm{ext}}\, S_{k_{\mathrm{ext}}}\tau$ in the exterior and $u_{\mathrm{int}} = D_{k_{\mathrm{int}}}\sigma + k_{\mathrm{ext}}\, S_{k_{\mathrm{int}}}\tau$ in the interior, the transmission conditions $u_{\mathrm{ext}} - u_{\mathrm{int}} = g$ and $\partial u_{\mathrm{ext}}/\partial\bm{n} - \partial u_{\mathrm{int}}/\partial\bm{n} = k_{\mathrm{ext}}\, h$ on $S$ lead to
\begin{equation}
\begin{aligned}
	\sigma + (D_{k_{\mathrm{ext}}} - D_{k_{\mathrm{int}}})\sigma + k_{\mathrm{ext}}\,(S_{k_{\mathrm{ext}}} - S_{k_{\mathrm{int}}})\tau &= g\,,\\
	-\tau + \tfrac{1}{k_{\mathrm{ext}}}(T_{k_{\mathrm{ext}}} - T_{k_{\mathrm{int}}})\sigma + (D'_{k_{\mathrm{ext}}} - D'_{k_{\mathrm{int}}})\tau &= h\,,
\end{aligned}
\label{eq:muller}
\end{equation}
where the factor $k_{\mathrm{ext}}$ scales the second unknown. In modal form, writing $G^{\mathrm{e}}_m = G_{k_{\mathrm{ext}},m}$ and $G^{\mathrm{i}}_m = G_{k_{\mathrm{int}},m}$\,,
\begin{equation}
\begin{aligned}
	\sigma_m + 2\pi\int_{\gamma} \Big[ \frac{\partial (G^{\mathrm{e}}_m - G^{\mathrm{i}}_m)}{\partial \bm{n}'}\,\sigma_m + k_{\mathrm{ext}}\,(G^{\mathrm{e}}_m - G^{\mathrm{i}}_m)\,\tau_m \Big]\, r'\,\D\ell' &= g_m\,,\\
	-\tau_m + 2\pi\int_{\gamma} \Big[ \frac{1}{k_{\mathrm{ext}}}\frac{\partial^2 (G^{\mathrm{e}}_m - G^{\mathrm{i}}_m)}{\partial \bm{n}\,\partial \bm{n}'}\,\sigma_m + \frac{\partial (G^{\mathrm{e}}_m - G^{\mathrm{i}}_m)}{\partial \bm{n}}\,\tau_m \Big]\, r'\,\D\ell' &= h_m\,.
\end{aligned}
\label{eq:muller_modal}
\end{equation}%
The differencing reduces the hypersingular kernel to a log-singular one, so the same panel quadrature as in the CFIE case applies. Assembling the matrix therefore requires the full second-order derivatives from Section~\ref{sec:allderivatives}. The system is $2N \times 2N$ for each mode (two unknowns per point), and every entry is evaluated at both $k_{\mathrm{ext}}$ and $k_{\mathrm{int}}$\,, doubling the kernel evaluation cost per pair. The same derivative data are also required in axisymmetric Maxwell transmission formulations, such as \cite{lai2019fft}, so the present evaluator applies without modification.

We use manufactured solutions based on point sources to verify both solvers. For the CFIE, point sources placed inside the scatterer generate a known exterior field; the surface values of this field provide the Dirichlet boundary data, and the reconstructed exterior field is compared with the exact value at an exterior target. For M\"uller's formulation, point sources are placed both inside and outside the scatterer, with wavenumbers $k_{\mathrm{ext}}$ and $k_{\mathrm{int}}$\,, respectively. The corresponding jump data define a transmission problem, and the reconstructed fields are checked at exterior and interior targets. In both cases, the target fields are evaluated mode by mode using the same modal Green's function evaluator.

\begin{figure}[!ht]
\centering
\begin{subfigure}[b]{0.21\textwidth}
\includegraphics[width=\textwidth]{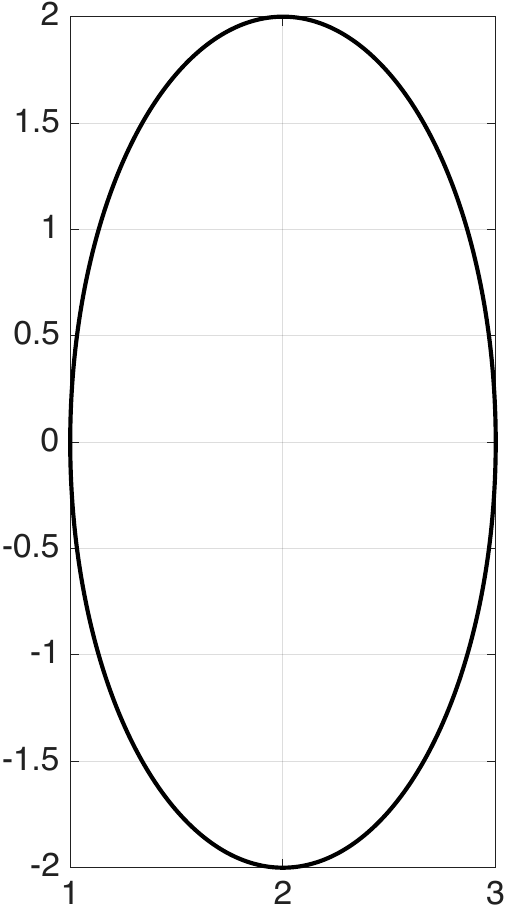}
\caption{}
\label{fig:torus_curve}
\end{subfigure}\qquad
\begin{subfigure}[b]{0.55\textwidth}
\includegraphics[width=\textwidth]{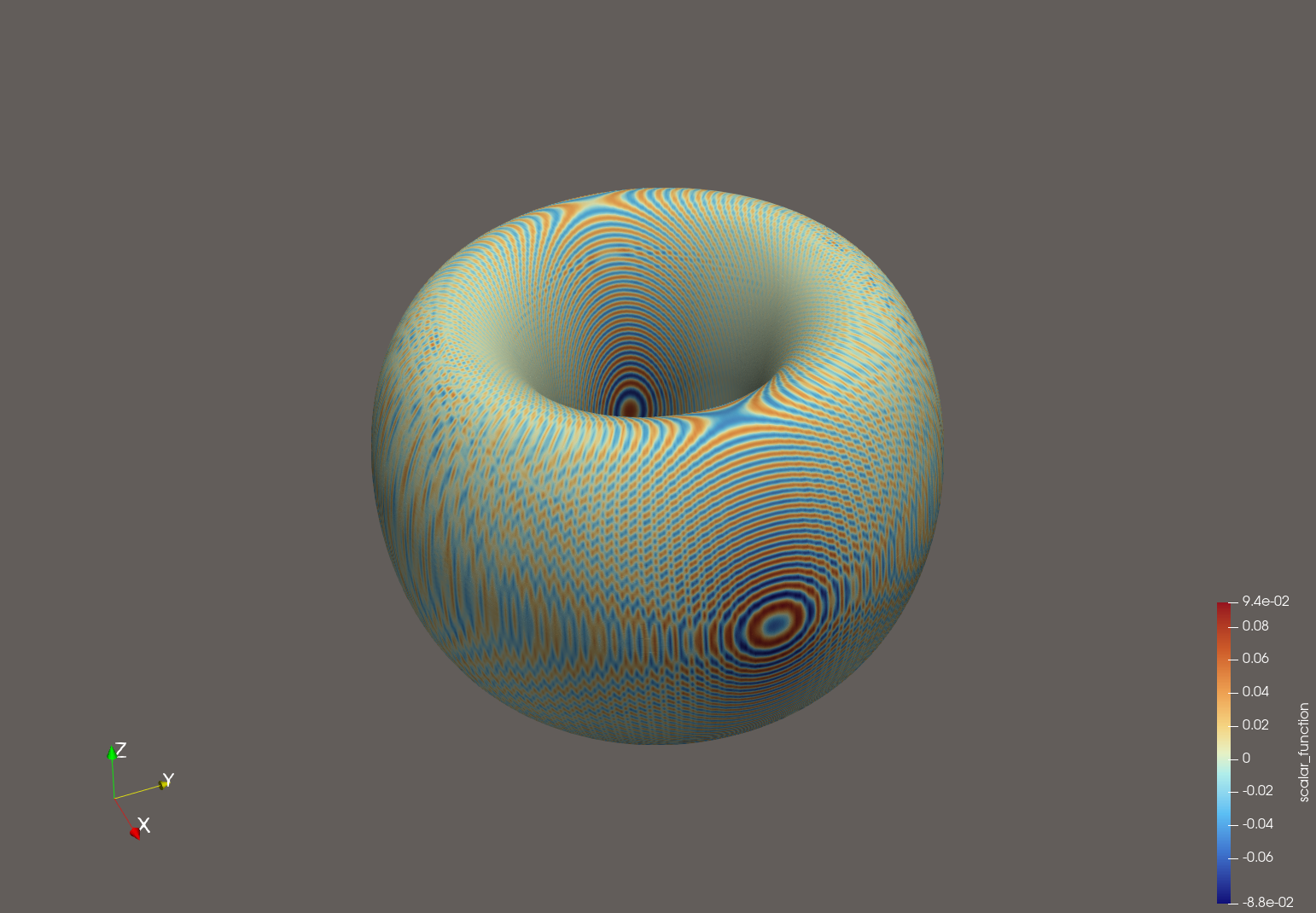}
\caption{}
\label{fig:torus_density}
\end{subfigure}
\caption{Sound-soft scattering from a torus ($k = 110$, $M = 330$). (a)~Generating curve in the $(r,z)$-plane, given by $r(t) = 2+\cos t$\,, $z(t) = 2\sin t$\,, $t \in [0, 2\pi)$\,. (b)~Real part of the density on the surface of revolution, with boundary data from two point sources inside the scatterer.}
\label{fig:torus}
\end{figure}

\begin{figure}[!ht]
\centering
\begin{subfigure}[b]{0.15\textwidth}
\includegraphics[width=\textwidth]{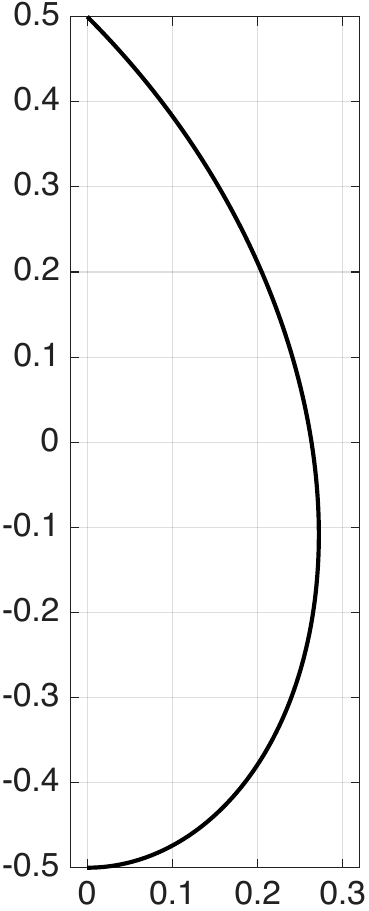}
\caption{}
\label{fig:droplet_curve}
\end{subfigure}\qquad
\begin{subfigure}[b]{0.55\textwidth}
\includegraphics[width=\textwidth]{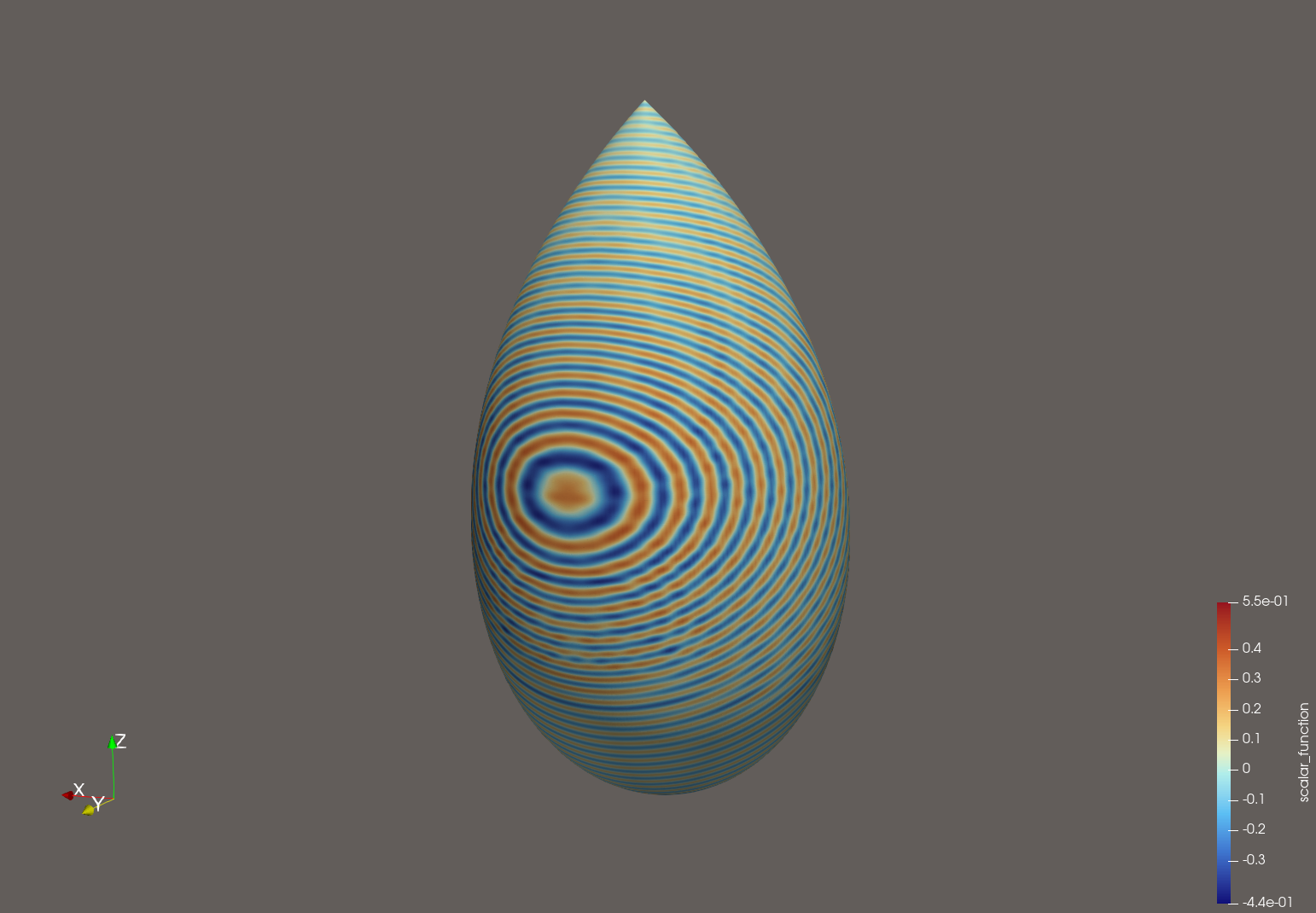}
\caption{}
\label{fig:droplet_density}
\end{subfigure}
\caption{Helmholtz transmission problem on a droplet ($k_{\mathrm{ext}} = 500$, $k_{\mathrm{int}} = 800$, $M = 150$). (a)~Generating curve in the $(r,z)$-plane, given by $r(t) = \cos(t/2)\sin(t/4)$\,, $z(t) = -\cos(t/2)\cos(t/4) + \tfrac{1}{2}$\,, $t \in [0, \pi]$\,. (b)~Real part of the double-layer density on the surface of revolution, with boundary data from point sources inside and outside the scatterer.}
\label{fig:droplet}
\end{figure}

The CFIE test uses a torus with elliptical cross-section (Figure~\ref{fig:torus}) at $k = 110$, while the transmission test uses a droplet-shaped body of revolution (Figure~\ref{fig:droplet}) at $k_{\mathrm{ext}} = 500$ and $k_{\mathrm{int}} = 800$\,. Both tests achieve more than $10$ digits of accuracy in the reconstructed field (Table~\ref{tab:bie}). The assembly is parallelized over source-target pairs and scales nearly linearly from one to eight cores; the LU solves are parallelized over Fourier modes.

\begin{table}[!ht]
\centering
\begin{subtable}{\textwidth}
\centering
\begin{tabular}{c|ccccccc}
\toprule
$k$ & $M$ & $N$ & Field error & \multicolumn{2}{c}{Assembly (s)} & LU solve (s) \\
    &     &     &             & 1-core & 8-core & \\
\midrule
110 & 330 & 2720 & $5.1 \times 10^{-12}$ & 545 & 91 & 56 \\
\bottomrule
\end{tabular}
\caption{CFIE on the torus of Figure~\ref{fig:torus}.}
\label{tab:bie_cfie}
\end{subtable}

\vspace{0.8em}

\begin{subtable}{\textwidth}
\centering
\begin{tabular}{cc|ccccccc}
\toprule
$k_{\mathrm{ext}}$ & $k_{\mathrm{int}}$ & $M$ & $N$ & Field error & \multicolumn{2}{c}{Assembly (s)} & LU solve (s) \\
 & & & & & 1-core & 8-core & \\
\midrule
500 & 800 & 150 & 1920 & $7.7 \times 10^{-11}$ & 458 & 76 & 78 \\
\bottomrule
\end{tabular}
\caption{M\"uller on the droplet of Figure~\ref{fig:droplet}.}
\label{tab:bie_muller}
\end{subtable}
\caption{Results for the BIE solvers. Assembly and LU solve times are totals across all Fourier modes $m = 0, \ldots, M$\,. The LU solves are parallelized over Fourier modes across 8 cores.}
\label{tab:bie}
\end{table}

\subsubsection{Wavenumber scaling}
\label{sec:bie_kscan}

The preceding tests validate the full modal BIE solve. We now isolate the asymptotic scaling of assembly and dense linear algebra by sweeping the wavenumber while keeping the number of Fourier modes fixed at $M = 40$\,. Since the evaluator has $k$-independent cost, increasing $k$ only increases the number $N$ of discretization points along the generating curve. Thus dense assembly should scale as $O(N^2)$ per mode, while the dense LU factorization should scale as $O(N^3)$\,.

For the CFIE sweep, we use the same torus geometry as in Figure~\ref{fig:torus}. The verification is performed mode by mode: the boundary data are generated by a modal point source, and the reconstructed modal field is checked at a target point. The reported field error is the maximum over all modes $m = 0, \ldots, M$ and is therefore larger than the single full-field error reported in Table~\ref{tab:bie}. The observed error floors in both sweeps are not primarily due to discretization: refining the CFIE case at $k = 100$ from $12$ to $20$ points per wavelength reduces the error only from $1.1 \times 10^{-9}$ to $3.3 \times 10^{-10}$\,, whereas a dominant 16th-order discretization error would decrease much more rapidly, and the modal systems are well conditioned (condition numbers about $30$). The level reflects the $11$--$12$ digit accuracy of the kernel evaluator (Table~\ref{tab:deriv_accuracy}). For M\"uller's formulation the condition numbers of the modal systems are larger, between $4 \times 10^{2}$ and $4 \times 10^{3}$ at $k_{\mathrm{ext}} = 500$\,, and the kernel differences in (\ref{eq:muller}) cost further digits through cancellation, which accounts for the lower accuracy in Table~\ref{tab:bie_kscan_muller}. The largest cases are limited by the laptop's $64$~GB of RAM rather than by runtime.

\begin{table}[!ht]
\centering
\begin{tabular}{r|rrrr}
\toprule
$k$ & $N$ & Field error & Assembly (s) & LU solve (s) \\
\midrule
100 & 1856 & $1.1 \times 10^{-9}$  &  24.7 &   2.3 \\
200 & 3712 & $3.0 \times 10^{-10}$ &  99.6 &  19.5 \\
300 & 5568 & $1.1 \times 10^{-9}$  & 224.5 &  72.7 \\
400 & 7408 & $7.1 \times 10^{-10}$ & 378.4 & 189.7 \\
\bottomrule
\end{tabular}
\caption{Wavenumber sweep for the CFIE solver on the geometry of Figure~\ref{fig:torus}, with $M = 40$ Fourier modes fixed. Assembly and LU solve times are totals across all modes, both run on 8 cores (assembly parallelized over source-target pairs, LU over Fourier modes). See Table~\ref{tab:bie} for the single-core assembly baseline.}
\label{tab:bie_kscan}
\end{table}

We repeat the experiment for M\"uller's formulation on the droplet geometry of Figure~\ref{fig:droplet}, using $k_{\mathrm{int}} = 1.6\,k_{\mathrm{ext}}$\,. The verification is again performed mode by mode, with modal point sources defining the transmission data and an exterior target used for the reported error.

\begin{table}[!ht]
\centering
\begin{tabular}{rr|rrrr}
\toprule
$k_{\mathrm{ext}}$ & $k_{\mathrm{int}}$ & $N$ & Field error & Assembly (s) & LU solve (s) \\
\midrule
 500 &  800 & 1920 & $1.6 \times 10^{-9}$ &  49.2 &  23.6 \\
 700 & 1120 & 2656 & $1.8 \times 10^{-8}$ &  97.0 &  61.9 \\
 900 & 1440 & 3408 & $6.6 \times 10^{-9}$ & 162.0 & 130.2 \\
1000 & 1600 & 3776 & $3.1 \times 10^{-9}$ & 199.3 & 195.3 \\
\bottomrule
\end{tabular}
\caption{Wavenumber sweep for the M\"uller solver on the droplet of Figure~\ref{fig:droplet}, with $M = 40$ and $k_{\mathrm{int}} = 1.6\,k_{\mathrm{ext}}$. The field error is the maximum over all modes of the relative reconstruction error at an exterior target. Assembly and LU solve times are totals across all modes; parallelization is as in Table~\ref{tab:bie_kscan}.}
\label{tab:bie_kscan_muller}
\end{table}

Tables~\ref{tab:bie_kscan}--\ref{tab:bie_kscan_muller} confirm the expected scaling. In both sweeps, the assembly time grows quadratically with $N$. This is a direct consequence of the $k$-independent kernel evaluator (Section~\ref{sec:accuracy_timing}); a kernel cost linear in $k$ would make dense assembly scale as $O(N^3)$\,. The dense LU cost grows as $O(N^3)$, as expected. For M\"uller on the droplet, the assembly-to-solve ratio decreases from $2.08$ to $1.02$, so the two costs are already comparable at $k_{\mathrm{ext}} = 1000$, beyond which the solve dominates. The CFIE sweep shows the same trend, with the ratio decreasing from $10.7$ to $2.0$, although the crossover lies beyond the range that fits in $64$~GB of RAM. Thus, in these experiments, kernel evaluation is no longer the limiting cost; the dominant remaining cost is dense per-mode linear algebra. Further acceleration should therefore come from applying fast direct solvers~\cite{martinsson2005fast,martinsson2019fast} based on hierarchical compression~\cite{cheng2005compression,ho2012fast} mode by mode.

\section{Conclusion}
\label{sec:conclusion}

In this paper, we have introduced an $O(M)$ algorithm for evaluating the azimuthal Fourier modes $G_m$\,, $m = 0, \ldots, M$\,, of the three-dimensional Helmholtz Green's function, together with all first- and second-order derivatives with respect to the cylindrical source and target coordinates. The cost is independent of the wavenumber and of the source-target separation.

The method combines three ingredients. First, we use the contour deformation of \cite{garritano2022efficient,gustafsson2010accurate} to evaluate a fixed number of boundary modes. The same deformed contour applies to the derivative integrals, and Generalized Gaussian Quadratures are used to handle the stronger near-singular factors that arise in the derivative integrands. Second, we cast the five-term recurrence satisfied by $G_m$ as a boundary-value problem in the mode index, rather than an unstable forward or backward recursion. This gives all intermediate modes in $O(M)$ time. A rigorous analysis of the componentwise accuracy observed across our experiments remains an interesting direction for future work. Third, the derivative data are obtained from stable recurrences, with cancellation-free combinations used in the near-singular regime.

The numerical experiments show that the resulting evaluator achieves high relative accuracy, typically close to double precision, for $G_m$\,, its first derivatives, and its second derivatives over a wide range of frequencies, separations, and mode numbers. In particular, relative accuracy is retained in the decay regime, where the modes are many orders of magnitude smaller than the integrand in the original oscillatory representation. The timing results confirm the expected linear scaling in $M$\,, and show that the runtime is essentially independent of both the wavenumber and the source-target separation.

We also demonstrated the use of the evaluator in modal boundary integral equation solvers for acoustic scattering from bodies of revolution. The examples include a combined-field integral equation for the exterior Dirichlet problem and M\"uller's formulation for transmission problems, the latter requiring second-order derivatives of the modal Green's function. In both cases, the computed fields agree with manufactured solutions to more than ten digits. The wavenumber sweeps show that, with the $k$-independent kernel evaluator, dense modal assembly scales quadratically with the number of discretization points along the generating curve, while the dominant remaining cost is dense per-mode linear algebra.

The algorithm extends naturally to complex wavenumbers $k$ with $\Im k > 0$ (e.g., attenuating media). As discussed in Remark~\ref{rem:complexk}, the steepest descent contours $\gamma_1$ and $\gamma_2$ rotate in the complex plane, with the new intersection points with the Bernstein ellipse obtained by Newton iteration at little additional cost. The recurrences of Section~\ref{sec:aux} and the pentadiagonal formulation of Section~\ref{sec:allmodes} are algebraic in $\kappa$ and apply to complex $\kappa$ as written; the remaining ingredient is the decay-regime threshold $m^*$ of (\ref{eq:mstar}) and the choice of $M'$ in Section~\ref{sec:decay}, which are derived for real $\kappa$ and need to be re-derived for complex $\kappa$\,, the Miller-type strategy itself being unchanged.

More broadly, with kernel evaluation no longer the bottleneck for axisymmetric BIE solvers, the natural next step is to combine the present evaluator with accelerated linear algebra, such as fast direct solvers based on hierarchical compression, applied mode by mode. Such a coupling should enable substantially larger frequency-domain scattering calculations while retaining the accuracy and geometric flexibility of boundary integral formulations.

\section*{Acknowledgments}
The author thanks Vladimir Rokhlin and Michael O'Neil for helpful discussions, and the anonymous referees for their helpful suggestions.

\newpage
\bibliographystyle{plain}
\bibliography{refs}

@article{gautschi1967computational,
	author = {Gautschi, Walter},
	journal = {SIAM Review},
	number = {1},
	pages = {24--82},
	title = {Computational aspects of three-term recurrence relations},
	volume = {9},
	year = {1967}}

@article{olver1967numerical,
	author = {Olver, F. W. J.},
	journal = {Journal of Research of the National Bureau of Standards, Section B: Mathematics and Mathematical Physics},
	number = {2--3},
	pages = {111--129},
	title = {Numerical solution of second-order linear difference equations},
	volume = {71B},
	year = {1967}}

@article{epstein2019high,
	author = {Epstein, Charles L. and Greengard, Leslie and O'Neil, Michael},
	journal = {Journal of Computational Physics},
	pages = {205--229},
	publisher = {Elsevier},
	title = {A high-order wideband direct solver for electromagnetic scattering from bodies of revolution},
	volume = {387},
	year = {2019}}

@article{matviyenko1995azimuthal,
	author = {Matviyenko, Gregory},
	journal = {Journal of Mathematical Physics},
	number = {9},
	pages = {5159--5169},
	publisher = {American Institute of Physics},
	title = {On the azimuthal {F}ourier components of the {G}reen's function for the {H}elmholtz equation in three dimensions},
	volume = {36},
	year = {1995}}

@techreport{gustafsson2010accurate,
	author = {Gustafsson, Mats},
	institution = {Lund University},
	number = {TEAT-7187},
	title = {Accurate and efficient evaluation of modal {G}reen's functions},
	year = {2010}}

@article{hao2014high,
	author = {Hao, S. and Barnett, A. H. and Martinsson, P. G. and Young, P.},
	journal = {Advances in Computational Mathematics},
	pages = {245--272},
	title = {High-order accurate methods for {N}ystr\"{o}m discretization of integral equations on smooth curves in the plane},
	volume = {40},
	year = {2014}}

@article{vico2014boundary,
	author = {Vico, Felipe and Greengard, Leslie and Gimbutas, Zydrunas},
	journal = {Numerische Mathematik},
	pages = {463--487},
	title = {Boundary integral equation analysis on the sphere},
	volume = {128},
	year = {2014}}

@article{lai2019fft,
	author = {Lai, Jun and O'Neil, Michael},
	journal = {Journal of Computational Physics},
	title = {An {FFT}-accelerated direct solver for electromagnetic scattering from penetrable axisymmetric objects},
	volume = {390},
	pages = {152--174},
	year = {2019},
	doi = {10.1016/j.jcp.2019.04.005}}

@article{kolm2001numerical,
	author = {Kolm, P. and Rokhlin, V.},
	journal = {Computers and Mathematics with Applications},
	number = {3-4},
	pages = {327--352},
	title = {Numerical quadratures for singular and hypersingular integrals},
	volume = {41},
	year = {2001}}

@article{bremer2010nonlinear,
	author = {Bremer, James and Gimbutas, Zydrunas and Rokhlin, Vladimir},
	journal = {SIAM Journal on Scientific Computing},
	number = {4},
	pages = {1761--1788},
	title = {A nonlinear optimization procedure for generalized {G}aussian quadratures},
	volume = {32},
	year = {2010},
	doi = {10.1137/080737046}}

@article{ho2012fast,
	author = {Ho, Kenneth L. and Greengard, Leslie},
	journal = {SIAM Journal on Scientific Computing},
	number = {5},
	pages = {A2507--A2532},
	publisher = {SIAM},
	title = {A fast direct solver for structured linear systems by recursive skeletonization},
	volume = {34},
	year = {2012}}

@article{cheng2005compression,
	author = {Cheng, Hongwei and Gimbutas, Zydrunas and Martinsson, Per-Gunnar and Rokhlin, Vladimir},
	journal = {SIAM Journal on Scientific Computing},
	number = {4},
	pages = {1389--1404},
	publisher = {SIAM},
	title = {On the compression of low rank matrices},
	volume = {26},
	year = {2005}}

@article{martinsson2005fast,
	author = {Martinsson, Per-Gunnar and Rokhlin, Vladimir},
	journal = {Journal of Computational Physics},
	number = {1},
	pages = {1--23},
	title = {A fast direct solver for boundary integral equations in two dimensions},
	volume = {205},
	year = {2005},
	doi = {10.1016/j.jcp.2004.10.033}}

@article{xue2023fullwave,
	author = {Xue, Wenjin and Zhang, Hanwen and Gopal, Abinand and Rokhlin, Vladimir and Miller, Owen D.},
	journal = {arXiv preprint arXiv:2308.08569},
	title = {Fullwave design of cm-scale cylindrical metasurfaces via fast direct solvers},
	year = {2023}}

@book{martinsson2019fast,
	author = {Martinsson, Per-Gunnar},
	title = {Fast Direct Solvers for Elliptic {PDE}s},
	publisher = {SIAM},
	series = {CBMS-NSF Regional Conference Series in Applied Mathematics},
	volume = {96},
	year = {2019},
	doi = {10.1137/1.9781611976045}}

@article{kressroach1978transmission,
	author = {Kress, R. and Roach, G. F.},
	title = {Transmission problems for the {H}elmholtz equation},
	journal = {Journal of Mathematical Physics},
	volume = {19},
	number = {6},
	pages = {1433--1437},
	year = {1978},
	doi = {10.1063/1.523808}}

@book{coltonkress2019inverse,
	author = {Colton, David and Kress, Rainer},
	title = {Inverse Acoustic and Electromagnetic Scattering Theory},
	publisher = {Springer},
	edition = {4th},
	year = {2019},
	series = {Applied Mathematical Sciences},
	volume = {93}}

@book{kress2014linear,
	author = {Kress, Rainer},
	title = {Linear Integral Equations},
	publisher = {Springer},
	edition = {3rd},
	year = {2014},
	series = {Applied Mathematical Sciences},
	volume = {82}}

@book{golub2013matrix,
	author = {Golub, Gene H. and Van Loan, Charles F.},
	title = {Matrix Computations},
	publisher = {Johns Hopkins University Press},
	edition = {4th},
	year = {2013}}

@article{chung2023inverse,
	author = {Chung, Haejun and Zhang, Feng and Li, Hao and Miller, Owen D. and Smith, Henry I.},
	journal = {Nanophotonics},
	number = {13},
	pages = {2371--2381},
	title = {Inverse design of high-{NA} metalens for maskless lithography},
	volume = {12},
	year = {2023}}

@article{oneil2018taylor,
	author = {O'Neil, Michael and Cerfon, Antoine J.},
	journal = {Journal of Computational Physics},
	pages = {263--282},
	title = {An integral equation-based numerical solver for {T}aylor states in toroidal geometries},
	volume = {359},
	year = {2018}}

@article{liu2016efficient,
	author = {Liu, Yu and Barnett, Alex H.},
	journal = {Journal of Computational Physics},
	pages = {226--245},
	title = {Efficient numerical solution of acoustic scattering from doubly-periodic arrays of axisymmetric objects},
	volume = {324},
	year = {2016}}

@article{youssef1989radar,
	author = {Youssef, Nashaat N.},
	journal = {Proceedings of the IEEE},
	number = {5},
	pages = {722--734},
	title = {Radar cross section of complex targets},
	volume = {77},
	year = {1989}}

@article{helsing2014kernel,
	author = {Helsing, Johan and Karlsson, Anders},
	journal = {Journal of Computational Physics},
	pages = {686--703},
	title = {An explicit kernel-split panel-based {N}ystr\"om scheme for integral equations on axially symmetric surfaces},
	volume = {272},
	year = {2014}}

@article{helsing2017resonances,
	author = {Helsing, Johan and Karlsson, Anders},
	journal = {IEEE Transactions on Microwave Theory and Techniques},
	number = {7},
	pages = {2214--2227},
	title = {Resonances in axially symmetric dielectric objects},
	volume = {65},
	year = {2017}}

@article{garritano2022efficient,
	author = {Garritano, James and Kluger, Yuval and Rokhlin, Vladimir and Serkh, Kirill},
	journal = {Journal of Computational Physics},
	pages = {111585},
	publisher = {Elsevier},
	title = {On the efficient evaluation of the azimuthal {F}ourier components of the {G}reen's function for {H}elmholtz's equation in cylindrical coordinates},
	volume = {471},
	year = {2022}}

@article{conway2010exact,
	author = {Conway, John T. and Cohl, Howard S.},
	title = {Exact {F}ourier expansion in cylindrical coordinates for the three-dimensional {H}elmholtz {G}reen function},
	journal = {Zeitschrift f{\"u}r angewandte Mathematik und Physik},
	volume = {61},
	number = {3},
	pages = {425--443},
	year = {2010},
	doi = {10.1007/s00033-009-0039-6}
}

\newpage
\appendix
\section{Numerically stable derivative formulas}
\label{sec:appendix_derivatives}

We collect here the chain rule formulas expressing all first- and second-order derivatives of $G_m$ with respect to the cylindrical coordinates $r,r',z,z'$ in terms of the $(a,b)$ derivatives, written in numerically stable form. We use the shorthand
\begin{align}
	&\Delta r = r - r'\,, \quad \Delta z = z - z'\,, \notag\\
	&S_1 = \frac{\partial^2 G_m}{\partial a^2} + \frac{\partial^2 G_m}{\partial a\partial b}\,, \quad S_2 = \frac{\partial^2 G_m}{\partial a\partial b} + \frac{\partial^2 G_m}{\partial b^2}\,.
\end{align}

\paragraph{First-order derivatives.}
\begin{align}
	\frac{\partial G_m}{\partial r} &= 2\Delta r\, \frac{\partial G_m}{\partial a} + 2r'\bb{\frac{\partial G_m}{\partial a} + \frac{\partial G_m}{\partial b}} \\
	\frac{\partial G_m}{\partial r'} &= -2\Delta r\, \frac{\partial G_m}{\partial a} + 2r\bb{\frac{\partial G_m}{\partial a} + \frac{\partial G_m}{\partial b}} \\
	\frac{\partial G_m}{\partial z} &= 2\Delta z\, \frac{\partial G_m}{\partial a} \\
	\frac{\partial G_m}{\partial z'} &= -\frac{\partial G_m}{\partial z}\,.
\end{align}

\paragraph{Second-order derivatives.}
\begin{align}
	\frac{\partial^2 G_m}{\partial r^2} &= 4r^2\, S_1 + 4r'^2\, S_2 - 4(\Delta r)^2\frac{\partial^2 G_m}{\partial a\partial b} + 2\frac{\partial G_m}{\partial a} \\
	\frac{\partial^2 G_m}{\partial r'^2} &= 4r'^2\, S_1 + 4r^2\, S_2 - 4(\Delta r)^2\frac{\partial^2 G_m}{\partial a\partial b} + 2\frac{\partial G_m}{\partial a} \\
	\frac{\partial^2 G_m}{\partial r\partial r'} &= 4rr'(S_1 + S_2) + 4(\Delta r)^2\frac{\partial^2 G_m}{\partial a\partial b} + 2\frac{\partial G_m}{\partial b} \\
	\frac{\partial^2 G_m}{\partial r\partial z} &= 4\Delta z\bb{r\, S_1 - \Delta r\,\frac{\partial^2 G_m}{\partial a\partial b}} \\
	\frac{\partial^2 G_m}{\partial r'\partial z} &= 4\Delta z\bb{r'\, S_1 + \Delta r\,\frac{\partial^2 G_m}{\partial a\partial b}} \\
	\frac{\partial^2 G_m}{\partial z^2} &= 2\frac{\partial G_m}{\partial a} + 4(\Delta z)^2\frac{\partial^2 G_m}{\partial a^2}
\end{align}
The remaining mixed derivatives are given by
\begin{align}
	\frac{\partial^2 G_m}{\partial r\partial z'} &= -\frac{\partial^2 G_m}{\partial r\partial z}\,, \\
	\frac{\partial^2 G_m}{\partial r'\partial z'} &= -\frac{\partial^2 G_m}{\partial r'\partial z}\,, \\
	\frac{\partial^2 G_m}{\partial z\partial z'} &= -\frac{\partial^2 G_m}{\partial z^2}\,, \\
	\frac{\partial^2 G_m}{\partial z'^2} &= \frac{\partial^2 G_m}{\partial z^2}\,.
\end{align}

\section{Contour integrals for derivative quantities}
\label{sec:appendix_contour}

We collect here the explicit contour integral decompositions for the derivative quantities of Section~\ref{sec:contour_integrands}. The notation follows that section: $u_- = \tau^2 - \I\beta_-$\,, $u_+ = \tau^2 - \I\beta_+$\,, and $z(\theta) = E_\rho(\theta)$ on the ellipse with $z'(\theta) = -a_e\sin\theta - \I b_e\cos\theta$\,.

\paragraph{First derivative $\frac{\partial G_m}{\partial a}$\,.}
From (\ref{eq:dgda_cheb}),
\begin{align}
	\frac{\partial G_m}{\partial a} = \frac{1}{8\pi^2 R_0^3}\bb{-I_{\gamma_1}^{(a)} + I_{\gamma_2}^{(a)} - I_{E_c}^{(a)}}\,,
	\label{eq:dgda_decomp}
\end{align}
where on $\gamma_1$\,:
\begin{align}
	I_{\gamma_1}^{(a)} = -\frac{4}{\sqrt{\alpha}} \int_{0}^{\tau_1} \frac{e^{-\kappa\sqrt{\alpha}\, u_-}}{\sqrt{\tau^2 - 2\I\beta_-}\,\sqrt{1+\gamma_1(\tau)}} \bb{\frac{\kappa}{\sqrt{\alpha}\,u_-} + \frac{1}{\alpha\, u_-^2}}\, T_m\bb{\gamma_1(\tau)}\, \D\tau\,,
	\label{eq:igamma1_dgda}
\end{align}
on $\gamma_2$\,:
\begin{align}
	I_{\gamma_2}^{(a)} = -\frac{4\I}{\sqrt{\alpha}} \int_{0}^{\tau_2} \frac{e^{-\kappa\sqrt{\alpha}\, u_+}}{\sqrt{\tau^2 - 2\I\beta_+}\,\sqrt{1-\gamma_2(\tau)}} \bb{\frac{\kappa}{\sqrt{\alpha}\,u_+} + \frac{1}{\alpha\, u_+^2}}\, T_m\bb{\gamma_2(\tau)}\, \D\tau\,,
	\label{eq:igamma2_dgda}
\end{align}
and on $E_c$\,:
\begin{align}
	I_{E_c}^{(a)} = \int_{\theta_1}^{\theta_2} \frac{e^{\I\kappa\sqrt{1-\alpha z}}}{\sqrt{1-z^2}} \bb{\frac{\I\kappa}{(1-\alpha z)} - \frac{1}{(1-\alpha z)^{3/2}}}\, T_m(z)\, z'(\theta)\,\D\theta\,.
	\label{eq:iec_dgda}
\end{align}

\paragraph{Combined first derivative $\frac{\partial G_m}{\partial a} + \frac{\partial G_m}{\partial b}$\,.}
From (\ref{eq:dgdapb_cheb}),
\begin{align}
	\frac{\partial G_m}{\partial a} + \frac{\partial G_m}{\partial b} = \frac{1}{8\pi^2 R_0^3}\bb{-I_{\gamma_1}^{(a+b)} + I_{\gamma_2}^{(a+b)} - I_{E_c}^{(a+b)}}\,,
	\label{eq:dgdapb_decomp}
\end{align}
where on $\gamma_1$\,:
\begin{align}
	I_{\gamma_1}^{(a+b)} = -\frac{4}{\sqrt{\alpha}} \int_{0}^{\tau_1} \frac{e^{-\kappa\sqrt{\alpha}\, u_-}}{\sqrt{\tau^2 - 2\I\beta_-}\,\sqrt{1+\gamma_1(\tau)}} \bb{\frac{\kappa}{\sqrt{\alpha}\,u_-} + \frac{1}{\alpha\, u_-^2}} (1-\gamma_1(\tau))\, T_m\bb{\gamma_1(\tau)}\, \D\tau\,,
\end{align}
on $\gamma_2$\,:
\begin{align}
	I_{\gamma_2}^{(a+b)} = -\frac{4\I}{\sqrt{\alpha}} \int_{0}^{\tau_2} \frac{e^{-\kappa\sqrt{\alpha}\, u_+}}{\sqrt{\tau^2 - 2\I\beta_+}\,\sqrt{1-\gamma_2(\tau)}} \bb{\frac{\kappa}{\sqrt{\alpha}\,u_+} + \frac{1}{\alpha\, u_+^2}} (1-\gamma_2(\tau))\, T_m\bb{\gamma_2(\tau)}\, \D\tau\,,
\end{align}
and on $E_c$\,:
\begin{align}
	I_{E_c}^{(a+b)} = \int_{\theta_1}^{\theta_2} \frac{e^{\I\kappa\sqrt{1-\alpha z}}}{\sqrt{1-z^2}} \bb{\frac{\I\kappa}{(1-\alpha z)} - \frac{1}{(1-\alpha z)^{3/2}}} (1-z)\, T_m(z)\, z'(\theta)\,\D\theta\,.
\end{align}
Here $(1-\gamma_1(\tau)) = -\tau^2(\tau^2 - 2\I\beta_-)$ and $(1-\gamma_2(\tau)) = 2 - \tau^4 + 2\I\beta_+\tau^2$\,.

\paragraph{Second derivative $\frac{\partial^2 G_m}{\partial a^2}$\,.}
From (\ref{eq:d2gda2_cheb}),
\begin{align}
	\frac{\partial^2 G_m}{\partial a^2} = \frac{1}{16\pi^2 R_0^5}\bb{-I_{\gamma_1}^{(aa)} + I_{\gamma_2}^{(aa)} - I_{E_c}^{(aa)}}\,.
	\label{eq:d2gda2_decomp}
\end{align}
The integrals have the same base as (\ref{eq:igamma1_gm})--(\ref{eq:iec_gm}), with the derivative factor $1/(\sqrt{1-\alpha z})$ replaced as follows. On $\gamma_1$\,:
\begin{align}
	I_{\gamma_1}^{(aa)} = -\frac{4}{\sqrt{\alpha}} \int_{0}^{\tau_1} \frac{e^{-\kappa\sqrt{\alpha}\, u_-}}{\sqrt{\tau^2 - 2\I\beta_-}\,\sqrt{1+\gamma_1(\tau)}} \bb{\frac{\kappa^2}{\alpha\, u_-^2} + \frac{3\kappa}{\alpha^{3/2}\, u_-^3} + \frac{3}{\alpha^2\, u_-^4}}\, T_m\bb{\gamma_1(\tau)}\, \D\tau\,.
\end{align}
On $\gamma_2$\,:
\begin{align}
	I_{\gamma_2}^{(aa)} = -\frac{4\I}{\sqrt{\alpha}} \int_{0}^{\tau_2} \frac{e^{-\kappa\sqrt{\alpha}\, u_+}}{\sqrt{\tau^2 - 2\I\beta_+}\,\sqrt{1-\gamma_2(\tau)}} \bb{\frac{\kappa^2}{\alpha\, u_+^2} + \frac{3\kappa}{\alpha^{3/2}\, u_+^3} + \frac{3}{\alpha^2\, u_+^4}}\, T_m\bb{\gamma_2(\tau)}\, \D\tau\,.
\end{align}
On $E_c$\,:
\begin{align}
	I_{E_c}^{(aa)} = \int_{\theta_1}^{\theta_2} \frac{e^{\I\kappa\sqrt{1-\alpha z}}}{\sqrt{1-z^2}} \bb{-\frac{\kappa^2}{(1-\alpha z)^{3/2}} - \frac{3\I\kappa}{(1-\alpha z)^2} + \frac{3}{(1-\alpha z)^{5/2}}}\, T_m(z)\, z'(\theta)\,\D\theta\,.
\end{align}

\paragraph{Combined second derivative $\frac{\partial^2 G_m}{\partial a^2} + \frac{\partial^2 G_m}{\partial a\partial b}$\,.}
From (\ref{eq:d2gda2pab_cheb}),
\begin{align}
	\frac{\partial^2 G_m}{\partial a^2} + \frac{\partial^2 G_m}{\partial a\partial b} = \frac{1}{16\pi^2 R_0^5}\bb{-I_{\gamma_1}^{(aa+ab)} + I_{\gamma_2}^{(aa+ab)} - I_{E_c}^{(aa+ab)}}\,,
\end{align}
where on $\gamma_1$\,:
\begin{align}
	I_{\gamma_1}^{(aa+ab)} &= -\frac{4}{\sqrt{\alpha}} \int_{0}^{\tau_1} \frac{e^{-\kappa\sqrt{\alpha}\, u_-}}{\sqrt{\tau^2 - 2\I\beta_-}\,\sqrt{1+\gamma_1(\tau)}} \notag\\
	&\qquad \cdot \bb{\frac{\kappa^2}{\alpha\, u_-^2} + \frac{3\kappa}{\alpha^{3/2}\, u_-^3} + \frac{3}{\alpha^2\, u_-^4}} (1-\gamma_1(\tau))\, T_m\bb{\gamma_1(\tau)}\, \D\tau\,,
\end{align}
on $\gamma_2$\,:
\begin{align}
	I_{\gamma_2}^{(aa+ab)} &= -\frac{4\I}{\sqrt{\alpha}} \int_{0}^{\tau_2} \frac{e^{-\kappa\sqrt{\alpha}\, u_+}}{\sqrt{\tau^2 - 2\I\beta_+}\,\sqrt{1-\gamma_2(\tau)}} \notag\\
	&\qquad \cdot \bb{\frac{\kappa^2}{\alpha\, u_+^2} + \frac{3\kappa}{\alpha^{3/2}\, u_+^3} + \frac{3}{\alpha^2\, u_+^4}} (1-\gamma_2(\tau))\, T_m\bb{\gamma_2(\tau)}\, \D\tau\,,
\end{align}
and on $E_c$\,:
\begin{align}
	I_{E_c}^{(aa+ab)} &= \int_{\theta_1}^{\theta_2} \frac{e^{\I\kappa\sqrt{1-\alpha z}}}{\sqrt{1-z^2}} \bb{-\frac{\kappa^2}{(1-\alpha z)^{3/2}} - \frac{3\I\kappa}{(1-\alpha z)^2} + \frac{3}{(1-\alpha z)^{5/2}}} \notag\\
	&\qquad \cdot (1-z)\, T_m(z)\, z'(\theta)\,\D\theta\,.
\end{align}

\end{document}